\documentclass{article}
\usepackage{graphicx}
\usepackage{psfrag}
\usepackage{amsfonts,amsmath,amssymb}
\usepackage{float}
\newtheorem{theorem}{Theorem}[section]
\newtheorem{lemma}[theorem]{Lemma}
\newtheorem{proposition}[theorem]{Proposition}
\newtheorem{corollary}[theorem]{Corollary}
\newtheorem{remark}[theorem]{Remark} 
\newenvironment{proof}{{\bf Proof\ }}{\QED\\}
\newcommand{\QED}{\hspace*{\fill}\rule{2.5mm}{2.5mm}}

\newcommand{\ds}{\displaystyle}

\newcommand{\R}{\mathbb{R}}
\newcommand{\N}{\mathbb{N}}

\newcommand{\ji}{J_i}
\newcommand{\jim}{J_{i-1}}
\newcommand{\nz} [1] {\vert\!\vert\!\vert\, #1 \vert\!\vert\!\vert}

\newcommand {\dtp}{D^{+}_{t}\,}
\newcommand {\dtm}{D^{-}_{t}\,}
\newcommand {\dtmm}{D^{--}_{t}\,}
\newcommand {\dtms}{D^{-*}_{t}\,}
\newcommand {\dtz}{D^{0}_{t}\,}
\newcommand {\dxp}{D^{+}_{x}\,}
\newcommand {\dxm}{D^{-}_{x}\,}

\newcommand {\dxz}{D^{0}_{x}\,}

\newcommand{\ie}[1]{[\negmedspace\lvert #1\rvert\negmedspace]}
\newcommand{\uu}{\mbox{{\footnotesize ${\cal U}$}}} 
\begin{document}

\title{Nonlinear Nonoverlapping Schwarz Waveform Relaxation for Semilinear Wave Propagation}

\author{Laurence Halpern\thanks{LAGA, Institut Galil{\'e}e,
  Universit\'e Paris XIII, 93430 Villetaneuse, France.
  halpern@math.univ-paris13.fr}
  \setcounter{footnote}{2}
  \and J\'er\'emie Szeftel\thanks{Department of Mathematics, Princeton University,
 Fine Hall, Washington Road, Princeton NJ 08544-1000, USA  and C.N.R.S., Math{\'e}matiques Appliqu{\'e}es de Bordeaux, Universit{\'e}  Bordeaux 1, 351 cours de la Lib{\'e}ration, 33405 Talence cedex France. jszeftel@math.princeton.edu.
The second author is partially supported by NSF Grant DMS-0504720}
}
\maketitle
\begin{abstract}
  We introduce a non-overlapping variant of the Schwarz waveform
  relaxation algorithm for semilinear wave propagation  in one
  dimension. Using the theory of absorbing boundary conditions, we
  derive a new nonlinear algorithm.  We show that the algorithm is
  well-posed and we prove its convergence by energy estimates and a
  Galerkin method. We then introduce an explicit scheme. We prove the
  convergence of the discrete algorithm with suitable assumptions on
  the nonlinearity. We finally illustrate our analysis with numerical
  experiments.
\end{abstract}
%\tableofcontents

% \begin{keywords}
%   Domain Decomposition, Waveform Relaxation, Schwarz Methods, Semilinear Wave Equation
% \end{keywords}

% \begin{AMS}
%   65F10, 65N22.
% \end{AMS}
%=================================================================
\section{Introduction}
\label{sec:intro}

%=================================================================
Schwarz waveform relaxation is a new class of algorithms for domain
decomposition in the frame of time dependant partial differential
equations. They are well-adapted to evolution problems, designed to
solve the equations separately on each spatial subdomain on the whole
time interval , exchanging informations on the space-time boundary of
the subdomains, overlapping or not \cite{Gander:1999:OCO}. In
particular for wave equations, it is of great importance, due to
numerical dispersion, to be able to handle local time and space
meshes, and this is allowed by the present method.  We presented the
method for the linear wave equation in \cite{GHN:2003:OSW} and
\cite{GH:2005:ABC} . When using overlapping subdomains and
``classical'' Schwarz waveform relaxation  -by a Dirichlet
  exchange of informations on the boundary- the so defined algorithm
converges in a finite number of iterations, inversely proportional to
the size of the overlap, which can be penalizing. We introduced
optimized transmission conditions, relying on the theory of absorbing
boundary conditions, which improve drastically the convergence of the
algorithm.

Very little has been done so far about nonlinear Schwarz algorithms.
An analysis of the classical Schwarz waveform relaxation algorithm was
performed in \cite{Gander:2003:OSW} for a conservation law.  The goal
of this paper is to define new Schwarz waveform relaxation algorithms
for the semilinear wave equation. We introduce two nonoverlapping
algorithms. The first one referred to as linear uses the absorbing
boundary condition of the linear problem whereas the second one
referred to as nonlinear uses the nonlinear absorbing boundary
conditions designed by J.
Szeftel in \cite{Sz:2006:ANL}. 

In Section 2, we introduce the definitions of the algorithms.

In Section 3 , we prove the algorithms to be well-posed. For the
precise analysis, we use a fixed point algorithm with regularity
estimates on a linear problem.

In Section 4, we prove the convergence of the algorithms. The proof is
an extension of a clever trick in \cite{Lions:1988:SAM}, already used
for linear algorithms, either hyperbolic or parabolic (see
\cite{GHN:2003:OSW}). However the nonlinearity requires a very fine
analysis.

In Section 5, we design discrete Schwarz waveform relaxation
algorithms. In each subdomain, the interior scheme is the usual
leapfrog scheme for the linear part, with a downwinding in time for
the nonlinear part.  The exchange of informations on the boundary is
naturally taken into account by a finite volume strategy.
In Section 6 we study the convergence of the algorithms, by discrete energy estimates.

As it is always the case for nonlinear problems, the well-posedness
and convergence results hold only locally in time. Therefore numerical
experiments are very important to bypass the limitations of the
theory. We present the results in Section 7, showing in particular
that our nonlinear algorithm gives optimal results within a large
class of algorithms.

\noindent\textbf{Remark}
  Due to the complexity of the mathematical theory, we restrain
  ourselves to the one dimensional case.  The multidimensional study
  contains additional difficulties due to the geometry and should be
  the heart of a forthcoming paper.

%=================================================================
\section{Problem Description}
\label{sec:pd}
%=================================================================

We consider the second order semilinear wave equation in one
dimension,
\begin{equation}\label{eq:WaveEquation}
(\partial^2_{t}-\partial^2_{x})\,\uu=f(\uu,\partial_t\, \uu,\partial_x\,\uu) 
\end{equation}
on the domain $\R\times(0,T)$ with initial conditions $\uu\,(\cdot,0)  =  p$, $\partial_t\,\uu\,(\cdot,0)  =  q$.\\

%
%--------------------------------------------------------------
\subsection{Absorbing Boundary Conditions for the Semilinear Wave Equation}
\label{subsec:abc}
%--------------------------------------------------------------

The question of absorbing boundary conditions arises when one wants to
make computations on an unbounded domain: 
a bounded computational domain is introduced, on
the boundary of which boundary conditions must be prescribed. These
boundary conditions must be absorbing to the waves leaving the domain.
A whole strategy has been designed by Engquist and Majda for linear
problems with variable coefficients, using pseudo-differential
operators \cite{EM:1979:RBC}. Recently it has been extended to
nonlinear operators by J.~Szeftel, in particular for the semilinear
wave equation \cite{Sz:2006:ANL}, using the paradifferential calculus
of \cite{Bo:1981:CSP} and \cite{Sa:1986:RMP}. We introduce a family
of operators
\begin{equation}
  \label{eq:transmoperator}
 {\cal B}^{\pm}(g^{\pm})u= \partial_tu \pm\partial_xu +g^{\pm}(u), 
\end{equation}
for   $C^{\infty}$ functions $g^{\pm}$ such that
$g^{\pm}(0)=0$. The linear absorbing boundary operators are given 
by $g^{\pm}=0$. In the case where $f(u,u_t,u_x)=f_1(u)+f_2(u)u_t+f_3(u)u_x$ with $f_j$ in $C^{\infty}(\R)$, $1\le j\le 3$, and $f_1(0)=0$, the following nonlinear boundary operators are given in \cite{Sz:2006:ANL}:
\begin{equation}\label{eq:para}
g^{+}(u):=-\frac{1}{2}\int_0^u(f_2-f_3)(\xi)d\xi, \qquad g^{-}(u):=-\frac{1}{2}\int_0^u(f_2+f_3)(\xi)d\xi.
\end{equation}
We replace the problem on the domain $\R$ by a boundary value problem in $\Omega_0=(a,b)$:
\begin{equation}\label{clanonlin} 
\begin{array}{c} 
(\partial^2_{t}-\partial^2_x)\bar{u}=f(\bar{u},\partial_t\bar{u},\partial_x\bar{u})  \mbox{ in $\Omega_0 \times (0,T)$}, \\[1mm] 
{\cal B}^{-}(g^-)\bar{u}(a,\cdot)=0,\quad{\cal B}^{+}(g^+)\bar{u}(b,\cdot)=0,
\end{array} 
\end{equation}
with initial values $p$ and $q$. Such boundary conditions give well-posed initial boundary value
problems, and are absorbing provided the intial data be compactly supported in $\Omega_0$, see
\cite{Sz:2006:ANL}. Following the strategy in \cite{GHN:2003:OSW}, we
use such absorbing operators for domain decomposition.

%------------------------------------------------------------
\subsection{A General Non-Overlapping Schwarz Waveform Relaxation Algorithm}
\label{subsec:gno}
%--------------------------------------------------------------
We decompose the domain $(a,b)$ into $I$ non overlapping subdomains
$\Omega_i=(a_i,a_{i+1})$, $a_j<a_i$ for $j<i$ and $a_1=a$,
$a_{I+1}=b$, and we introduce a general non overlapping Schwarz waveform relaxation algorithm. An initial guess $\{h^{\pm,0}_i\}_{1\le i \le I+1}$ is given. For $k \ge 1$, one step of the algorithm is 
\begin{equation}\label{eq:NonoverlappingSchwarz}
\begin{array}{l} 
\left[
\begin{array}{l} 
\ds   (\partial^2_{t}-\partial^2_x)u_i^{k} = f(u_i^{k},\partial_tu_i^{k},\partial_xu_i^{k}) \mbox{ in $\Omega_i \times (0,T)$},\\
\ds     u_i^{k}(\cdot,0) = p,\quad
    \partial_{t}u_i^{k}(\cdot,0) = q \mbox{ in }\Omega_i,\\  
\ds     {\cal B}^-(g^-)u^{k}_i(a_i,\cdot) = h^{-,\,k-1}_i,\quad
    {\cal B}^+(g^+)u^{k}_i(a_{i+1},\cdot)  =h^{+,\,k-1}_i \mbox{ in }(0,T) ,\\
\end{array}
\right.
 \\[8mm]
\ds    h^{-,\, k}_i= {\cal B}^-(g^-)u_{i-1}^{k}(a_i,\cdot),\quad
   h^{+,\, k}_i ={\cal B}^+(g^+)u_{i+1}^{k}(a_{i+1},\cdot) \mbox{ in }(0,T) ,\\
\end{array}
\end{equation}
where ${\cal B}^{\pm} $ are given in (\ref{eq:transmoperator}). For
ease of notations, we defined here $h_{1}^{\pm,\, k}=0$ and $h_{I+1}^{\pm,\, k}=0$, so
that the index $i$ in (\ref{eq:NonoverlappingSchwarz}) ranges from
$i=1,2,\ldots,I$. In 
the sequel, we call linear transmission condition the choice
$g^{\pm}=0$ and nonlinear transmission condition the choice
\eqref{eq:para}. For the
classical linear homogeneous wave equation, it has been proved in
\cite{GHN:2003:OSW} that the algorithm converges optimally if $T$ is
small enough (which means in two iterations, independently of the
number of subdomains), and the transmission operators ${\cal B}^\pm$
are given by ${\cal B}^\pm=\partial_t\pm
\partial_x$. This behavior is due to the finite speed of propagation, together with the fact that these operators are the exact Dirichlet Neumann operators in this case. In the nonlinear case, the propagation still takes place with the finite speed, but we can use only approximate Dirichlet Neumann operators. Therefore the classical Schwarz algorithm with overlap is still convergent, and for our nonoverlapping nonlinear algorithms, we will use energy estimates.

%=================================================================
\section{Well-posedness For The Subproblems}
\label{sec:wp}
%=================================================================
The study of the nonlinear problem relies on an iterative linear scheme. Therefore a first step for the definition of the algorithm is the study of the nonhomogeneous initial boundary value problem for a general domain $\Omega=(a_-,a_+)$,
\begin{equation}
\label{eq:IBVPNL}
\begin{array}{c} 
(\partial^2_{t}-\partial^2_x)u=f(u,\partial_tu,\partial_xu)\,\, \mbox{in $\Omega \times (0,T)$}, \\[1mm] 
{\cal B}^{-}(g^-)u(a_-,\cdot)=h^-,\quad
{\cal B}^{+}(g^+)u(a_+,\cdot)=h^+,
\end{array} 
\end{equation}
with initial values $p$ and $q$. We will use for $j \le 2$ the spaces
\begin{equation}
  \label{eq:defv_j}
V_j(\Omega,T)=\{u \in L^{\infty}(0,T;L^{2}(\Omega)),\partial^{\alpha}u \in  L^{\infty}(0,T;L^{2}(\Omega)), |\alpha|\le j\} . 
\end{equation}
In formula (\ref{eq:defv_j}), $\alpha$ is a 2-index in $\N^2$, the first coordinate in $\alpha$ stands for the time, and the second one stands for the space, so for instance $\partial^{\alpha}u=\partial_{tx}u$ for $\alpha=(1,1)$. $V_j(\Omega,T)$ is equipped with the norm $\|u\|_{V_j(\Omega,T)}= \max_{|\alpha|\le j}\|\partial^{\alpha}u\|_{L^{\infty}(0,T;L^{2}(\Omega))}$. 
\begin{theorem}\label{th:wpibvpnl}
  Let $p$ in $H^{2}(\Omega)$ and $q$ in $H^{1}(\Omega)$. There exists
  a time $T^*$ such that for any $T\le T^*$, for $h^\pm$ in $H^1(0,T)$ with the
  compatibility conditions 
\begin{equation}
  \label{eq:cc}
h^\pm(0)= q(a_\pm)\pm p'(a_\pm)+g^\pm(p(a_\pm)),
\end{equation}
(\ref{eq:IBVPNL}) has a unique solution $u$ in $V_2(\Omega,T)$, with
$\partial_tu(a_\pm,\cdot)$ and $\partial_xu(a_\pm,\cdot)$ in
$H^1(0,T)$. Furthermore there exists a positive real number $C^*$ such
that
\begin{equation}\label{eq:estnl}
\|u\|^2_{V_2(\Omega,T)}+\ds\sum_{\pm}\sum_{|\alpha|=1}\|\partial^\alpha u(a_\pm,\cdot)\|^2_{H^1(0,T)} 
\le C^*(\|p\|^2_{H^2(\Omega)}+\|q\|^2_{H^1(\Omega)}+ \ds\sum_{\pm}\|h^\pm\|^2_{H^1(0,T)}),
\end{equation}
where $T^*$ and  $C^*$ depend on the data $p,q,f,g^\pm,h^\pm$.
\end{theorem}
This result has been first proved in \cite{Sz:2006:ANL} with homogeneous boundary conditions ({\it i.e.} $h^\pm=0$). The additional difficulty  comes from the boundary conditions, and we give here the main steps of the proof. It relies on the construction of a sequence of linear problems of the form: 
\begin{equation}
\label{eq:IBVP}
  \begin{array}{c}
    \partial^2_{t}\tilde{u}-\partial^2_x\tilde{u}+\tilde{u} = F \,\, \mbox{in $\Omega \times (0,T)$},\\
    (\partial_t\tilde{u} -\partial_x\tilde{u} )(a_-,\cdot) = H^-,\quad
    (\partial_t\tilde{u} +\partial_x\tilde{u} )(a_+,\cdot) = H^+.
 \end{array}
\end{equation}
\begin{proposition}\label{prop:wpibvplin}
  Let $p$ in $H^{2}(\Omega)$ and $q$ in $H^{1}(\Omega)$. For any positive time $T$, let $F$ in
  $H^1((0,T)\times\Omega)$, and  $H^\pm$ in $H^1(0,T)$ with the
  compatibility conditions 
\begin{equation}\label{eq:ccl}
H^\pm(0)= q(a_\pm)\pm p'(a_\pm).
\end{equation}
Then, (\ref{eq:IBVP}) with initial data $p$ and $q$ has a unique
solution $\tilde{u}$ in $V_2(\Omega,T)$, with $\partial_t\tilde{u}(a_\pm,\cdot)$ and
$\partial_x\tilde{u}(a_\pm,\cdot)$ in $H^1(0,T)$.  Moreover we have the following bounds on the solution
\begin{multline}\label{eq:estimsecond}
\|\tilde{u}\|^2_{V_2(\Omega,T)}
+ \ds\sum_{\pm}\bigl(\|\partial_t\tilde{u}(a_\pm,\cdot)\|^2_{H^1(0,T)} +\|\partial_x\tilde{u}(a_\pm,\cdot)\|^2_{H^1(0,T)}\bigr)\\ 
 \le C_1e^T\bigl(\|F\|^2_{H^1(0,T;L^{2}(\Omega))}+\|F(\cdot,0)\|^2_{L^{2}(\Omega)}\
+\ds\sum_{\pm}\|H^\pm\|^2_{H^1(0,T)}+\|p\|^2_{H^2(\Omega)}+\|q\|^2_{H^1(\Omega)}\bigr),
\end{multline}
where $C_1$ is a universal constant.
\end{proposition}
\begin{proof}
  We start with the {\it a priori} estimates. We multiply
  (\ref{eq:IBVP}) by $\partial_t\tilde{u}$ and integrate by parts in $\Omega$:
\begin{multline*}
\ds\frac{1}{2}\ds\frac{d}{dt}\bigl(\|\partial_t\tilde{u}(\cdot,t)\|_{L^{2}(\Omega)}^2+ \|\tilde{u}(\cdot,t)\|^2_{H^{1}(\Omega)}\bigr)
+ (\partial_t\tilde{u}(a_-,t))^2+ (\partial_t\tilde{u}(a_+,t))^2 \\[3mm]
  = (F(\cdot,t), \partial_t\tilde{u}(\cdot,t))_{L^{2}(\Omega)}+ H^-(t)\partial_t\tilde{u}(a_-,t)+ H^+(t)\partial_t\tilde{u}(a_+,t). 
\end{multline*}
Using the Cauchy-Schwarz inequality on the right hand side, together
with the inequality $ \alpha\beta\le \frac{1}{2}\ \alpha^2+\frac{1}{2
}\ \beta^2$ for all $\alpha,\beta\in \R$, and finally integrating in
time, we obtain
\begin{multline*}
\ds\|\partial_t\tilde{u}(\cdot,t)\|_{L^{2}(\Omega)}^2+ \|\tilde{u}(\cdot,t)\|^2_{H^{1}(\Omega)}+\int_0^t[(\partial_t\tilde{u}(a_-,s))^2+ (\partial_t\tilde{u}(a_+,s))^2]\,ds \\
\le  \ds\int_0^t[\|\partial_t\tilde{u}(\cdot,s)\|_{L^{2}(\Omega)}^2\,ds+ \int_0^t\|F(\cdot,s)\|_{L^{2}(\Omega)}^2\,ds\\
+\|q\|_{L^{2}(\Omega)}^2+ \|p\|^2_{H^{1}(\Omega)}+ \ds\int_0^t[(H^-(s))^2 +(H^+(s))^2 ]\,ds.
\end{multline*}
By Gronwall Lemma, we deduce that 
\begin{multline*}
\ds\|\partial_t\tilde{u}(\cdot,t)\|_{L^{2}(\Omega)}^2+ \|\tilde{u}(\cdot,t)\|^2_{H^{1}(\Omega)}+\int_0^t[(\partial_t\tilde{u}(a_-,s))^2+ (\partial_t\tilde{u}(a_+,s))^2]\,ds \\
\le   \ds   e^T\bigl(\|F\|^2_{L^2((0,T)\times\Omega)}
+\|q\|_{L^{2}(\Omega)}^2+ \|p\|^2_{H^{1}(\Omega)}+ \ds\sum_{\pm}\|H^\pm\|_{L^{2}(0,T)}^2\bigr),
\end{multline*}
which gives
\begin{multline}
  \label{eq:FirstEnergyEstimate}
  \ds\max_{|\alpha|\le 1}\|\partial^{\alpha}\tilde{u}\|^2_{L^{\infty}(0,T;L^{2}(\Omega))}+\|\partial_t\tilde{u}(a_-,\cdot)\|_{L^2(0,T)}^2+ \|\partial_t\tilde{u}(a_+,\cdot)\|_{L^2(0,T)}^2 \\
  \le \ds e^T\bigl(\|F\|^2_{L^2((0,T)\times\Omega)}
  +\|q\|_{L^{2}(\Omega)}^2+ \|p\|^2_{H^1(\Omega)}+
  \ds\sum_{\pm}\|H^\pm\|_{L^2(0,T)}^2\bigr).
\end{multline}
Differentiating in time in (\ref{eq:IBVP}), we now apply
(\ref{eq:FirstEnergyEstimate}) to $\partial_t \tilde{u}$, and obtain:
\begin{multline}\label{eq:FirstbisEnergyEstimate}
\ds\max_{|\alpha|\le 1}\|\partial^{\alpha}\partial_t\tilde{u}\|^2_{L^{\infty}(0,T;L^{2}(\Omega))}+\|\partial_t^2\tilde{u}(a_-,\cdot)\|_{L^2(0,T)}^2+ \|\partial_t^2\tilde{u}(a_+,\cdot)\|_{L^2(0,T)}^2 \\
\le   \ds  e^T(\|\partial_tF\|^2_{L^2((0,T)\times\Omega)}
+\|\partial_t^2 u(\cdot,0)\|_{L^{2}(\Omega)}^2+ \|q\|^2_{H^1(\Omega)}+ \ds\sum_{\pm}\|\partial_tH^\pm\|_{L^2(0,T)}^2).
\end{multline}
We must estimate $\partial_t^2 \tilde{u}(\cdot,0)$ in the righthand side of \eqref{eq:FirstbisEnergyEstimate}. We multiply
(\ref{eq:IBVP}) by $\partial_t^2 \tilde{u}$, integrate in space and
evaluate at time 0:
\begin{multline*}
  \|\partial^2_t\tilde{u}(\cdot,0)\|^{2}_{L^{2}(\Omega)}+(\partial_xp,\partial_x\partial^2_t\tilde{u}(\cdot,0))+(p,\partial^2_t\tilde{u}(\cdot,0))\\
  +q(a_-)\partial^2_t\tilde{u}(a_-,0)+ q(a_+)\partial^2_t\tilde{u}(a_+,0) \\
  =H^-(0)\partial^2_t\tilde{u}(a_-,0)+H^+(0)\partial^2_t\tilde{u}(a_+,0)+(F(\cdot,0),\partial^2_t\tilde{u}(\cdot,0)).
\end{multline*}
We integrate by parts in the second term, and rewrite the equality as
\begin{multline*}
\|\partial^2_t\tilde{u}(\cdot,0)\|^{2}_{L^{2}(\Omega)}-(\partial_x^2\tilde{u}(\cdot,0),\partial^2_t\tilde{u}(\cdot,0))+(p,\partial^2_t\tilde{u}(\cdot,0))\\
=(H^-(0)-q(a_-)+\partial_xp(a_-))\partial^2_t\tilde{u}(a_-,0)+(H^+(0)-q(a_+)-\partial_xp(a_+))\partial^2_t\tilde{u}(a_+,0)+(F(\cdot,0),\partial^2_t\tilde{u}(\cdot,0)).
\end{multline*}
The boundary terms on the right-hand side vanish by the compatibility conditions, and we get
\[
\begin{array}{l}   
\|\partial^2_t\tilde{u}(\cdot,0)\|^{2}_{L^{2}(\Omega)}=(\partial_x^2p-p+F(\cdot,0),\partial^2_t\tilde{u}(\cdot,0)).
\end{array} 
\]
Using the Cauchy-Schwarz Lemma, we obtain
\[
\|\partial^2_t\tilde{u}(\cdot,0)\|_{L^{2}(\Omega)} \le \|\partial_x^2p-p+F(\cdot,0)\|_{L^{2}(\Omega)}.
\]
We replace the term $\|\partial^2_t\tilde{u}(\cdot,0)\|_{L^{2}(\Omega)}$ in \eqref{eq:FirstEnergyEstimate}, and we deduce the second a priori estimate: 
\begin{multline}
  \label{eq:SecondEnergyEstimate}
\ds\max_{|\alpha|\le 1}\|\partial^{\alpha}\partial_t\tilde{u}\|^2_{L^{\infty}(0,T;L^{2}(\Omega))}+\|\partial_t^2\tilde{u}(a_-,\cdot)\|_{L^2(0,T)}^2+ \|\partial_t^2\tilde{u}(a_+,\cdot)\|_{L^2(0,T)}^2 \\
\le   \ds 3\, e^T\bigl(\|\partial_tF\|^2_{L^2((0,T)\times \Omega)}
+\|F(\cdot,0)\|^2_{L^{2}(\Omega)}+\| p\|_{H^{2}(\Omega)}^2+ \|q\|^2_{H^{1}(\Omega)}+ \ds\sum_{\pm}\|\partial_tH^\pm\|_{L^2(0,T)}^2\bigr).
\end{multline}
We still need to estimate the mixed derivatives $\partial_{xx}\tilde{u}$ in
the interior and $\partial_{xt}\tilde{u}$ on the boundaries. We use the
equation, which gives in the interior
\[
\|\partial_{xx} \tilde{u}\|_{L^{\infty}(0,T;L^{2}(\Omega))}\le \|\partial_{tt} \tilde{u}\|_{L^{\infty}(0,T;L^{2}(\Omega))}+ \|\tilde{u}\|_{L^{\infty}(0,T;L^{2}(\Omega))}+ \|F\|_{L^{\infty}(0,T;L^{2}(\Omega))}.
\]
We now introduce the inequality
\[
\|F\|^2_{L^{\infty}(0,T;L^{2}(\Omega))} \le 2(\|F(\cdot,0)\|^2_{L^{2}(\Omega)}+ T\|\partial_t F\|^2_{L^{2}(0,T;L^{2}(\Omega))} ),
\]
and by (\ref{eq:SecondEnergyEstimate}) and (\ref{eq:FirstEnergyEstimate}), we get, using that $e^T \ge 1$ and $e^T \ge T$,
\begin{multline*}
\|\partial_{xx} \tilde{u}\|^2_{L^{\infty}(0,T;L^{2}(\Omega))}\le 15 e^T\bigl(\|F\|^2_{H^{1}(0,T;L^{2}(\Omega))} + \|F(\cdot,0)\|^2_{L^{2}(\Omega)}
            + \| p\|_{H^{2}(\Omega)}^2+ \|q\|^2_{H^{1}(\Omega)}  +\sum_{\pm}\| H^\pm\|_{H^{1}(0,T)}^2 \bigr).
\end{multline*}
As for the boundary term, we get for instance on the left boundary
\[
\|\partial_{tx}\tilde{u}(a_-,\cdot)\|_{L^2(0,T)}\le \|\partial_{tt}\tilde{u}(a_-,\cdot)\|_{L^2(0,T)}+ \|\partial_{t}H^-\|_{L^2(0,T)}.
\]
Squaring the inequality, and adding the term coming from the right boundary leads to
\[
\ds\sum_{\pm}\|\partial_{tx}\tilde{u}(a_\pm,\cdot)\|^2_{L^2(0,T)}\le
2\ds\sum_{\pm}\bigl(\|\partial_{tt}\tilde{u}(a_\pm,\cdot)\|^2_{L^2(0,T)}
+\| \partial_{t}H^\pm\|^2_{L^2(0,T)}\bigr),
\]
which provides the last estimate announced in the proposition. The well-posedness is then derived in a standard way  by the Galerkin method.
\end{proof}
The solution $\bar{u}$ of the nonlinear subdomain problem is now defined
through an iterative scheme. 
The initial guess is $\bar{u}_0=p$.
At step $k$, $\bar{u}_{k}$ being known, we define
\begin{equation}
  \label{eq:SecondMembrek}
\begin{array}{lcl}
{\cal F}(w)=w+f(w,\partial_t w, \partial_x w),\ 
{\cal G}^\pm (w) =g^\pm (w(a_\pm,\cdot)),\ 
{\cal H}^\pm(w) =h^\pm-{\cal G}^\pm (w).
\end{array}
\end{equation}
$\bar{u}_{k+1}$ is the solution of the linear initial boundary value problem (\ref{eq:IBVP}) with data $f_k={\cal F}(\bar{u}_k)$, $g_k^\pm={\cal G}^\pm (\bar{u}_k)$, $h_k^\pm={\cal H}(\bar{u}_k)$, and initial data $p$ and $q$.  The proof of convergence for the sequence $\bar{u}_k$ is written in details in \cite{Sz:2006:ANL}. The uniqueness follows from the result:
\begin{lemma}\label{lem:estimfG}
  There exists a real positive increasing function $\theta$ such that,
  for any time $T$, for any $v$  in
  $V_2(\Omega,T)$, $\partial^{\alpha}{\cal F}(v)$ is in
  $V_1(\Omega,T)$,  ${\cal G}^\pm (v)$ and ${\cal H}(v)$ are in
  $L^{\infty}(0,T)$. Moreover, for $v_1$, $v_2$ in  $V_2(\Omega,T)$, we have
\begin{equation}\label{eq:estimfG}
\begin{array}{c}
\ds\|{\cal H}^\pm (v_1)-{\cal H}^\pm (v_2)\|_{L^{\infty}(0,T)} \le \theta( \|v_1\|_{V_2(\Omega,T)}+\|v_2\|_{V_2(\Omega,T)}) \|v_1-v_2\|_{V_1(\Omega,T)},\\[3mm]
\ds\|{\cal F} (v_1)-{\cal F} (v_2)\|_{V_1(\Omega,T)} \le\theta( \|v_1\|_{V_2(\Omega,T)}+\|v_2\|_{V_2(\Omega,T)})\|v_1-v_2\|_{V_2(\Omega,T)},
\end{array}
\end{equation}
\end{lemma}

As a consequence, we have the well-posedness of problem (\ref{clanonlin}).
\begin{corollary}
  Let $p$ in $H^{2}_0(\Omega)$ and $q$ in $H^{1}_0(\Omega)$. There
  exists a time $T_0^*$ such that for any $T\le T_0^*$,
  (\ref{clanonlin}) has a unique solution $\bar{u}$ in
  $V_2(\Omega,T)$, with $\partial^\alpha\bar{u}(a,\cdot)$ and
  $\partial^\alpha\bar{u}(b,\cdot)$ in $H^1(0,T)$ for
  $\lvert\alpha\rvert=1$. Furthermore there exists a positive real
  number $C^*$ depending only on the size of $\Omega$ such that
\begin{equation*}\label{eq:estnltot}
\begin{array}{c}
\|\bar{u}\|^2_{V_2(\Omega,T)}+\ds\sum_{|\alpha|=1}\|\partial^\alpha \bar{u}(a,\cdot)\|^2_{H^1(0,T)}\ds\sum_{|\alpha|=1}\|\partial^\alpha \bar{u}(b,\cdot)\|^2_{H^1(0,T)}
\le C^*(\|p\|^2_{H^2(\Omega)}+\|q\|^2_{H^1(\Omega)}).
\end{array}
\end{equation*}
\end{corollary} 

%=================================================================
\section{Convergence of The Algorithm}
\label{sec:conv}
%=================================================================
%
We now study the convergence of the Schwarz waveform Relaxation Algorithm (\ref{eq:NonoverlappingSchwarz}). In order to define the algorithm, we need a regularity result:
\begin{proposition}\label{th:regularity}
For any $\epsilon$, $0 < \epsilon < 1$, $V_2(\Omega,T)\subset  {\cal C}^0(0,T;H^{2-\epsilon}(\Omega))\cap  {\cal C}^1(0,T;H^{1-\epsilon}(\Omega))$.
\end{proposition}
\begin{proof}
By using extension operators in time and space, it suffices to prove the result in $\R \times \R$. We make use of the Littlewood-Paley theory (see for example \cite{Che:1995:MF}). In particular, there exists $\varphi$ and $\chi$ two tempered distributions on $\R$, with $\varphi$ supported in $(-8/3,-3/4)\cup(3/4,8/3)$, $\chi$ supported in $(-4/3,4/3)$, and  
\[
\chi(\xi)+\sum_{q\ge 0}\varphi(2^{-q}\xi) =1,\,\forall\xi\in\R.
\]
We define the dyadic projectors $\Delta_q$ by their action on a function u,
\begin{equation}\label{eq:deltaq}
\Delta_{-1}u=\chi(D)u, \,\Delta_{q}u=\varphi(2^{-q}D)u \mbox{ for } q\ge 0,
\end{equation}
where $D=-i\partial$. These operators give an equivalent norm in $H^s(\R)$,
\[
\vert u \vert_s=\biggl(\sum_{q \ge -1}2^{2qs}\Vert\Delta_q u\Vert^2 \biggr)^{\frac{1}{2}}.
\]
They can also be used to define the Zygmund spaces
\[
{\cal C}^r_*=\{ u \in {\cal S}', \nz{u}_r=
\sup_{q \ge -1}2^{qr}\Vert\Delta_q u\Vert^2 < +\infty\}.
\]
${\cal C}^r_*$ co\"incides with the usual H\"older space when $r$ is not an
integer. For any positive $r$, we know that $W^{r,\infty}$, the space of
functions in $L^\infty$ with derivatives of order up to $r$ in
$L^\infty$, is included in ${\cal C}^r_*$. Therefore we have
\[
V_2(\R,\R)\subset {\cal C}^0_*(\R,H^2(\R))\cap  {\cal C}^1_*(\R,H^1(\R))\cap  {\cal C}^2_*(\R,L^2(\R)).
\]
We need an interpolation lemma.
\begin{lemma} \label{lem:interp}
For any positive $\alpha,\beta,a,b$, for any $\theta$, $0 \le \theta \le 1$, 
\[
{\cal C}^\alpha_*(\R,H^a(\R))\cap {\cal C}^\beta_*(\R,H^b(\R))\subset {\cal C}^{\theta\alpha+(1-\theta)\beta}_*(\R,H^{\theta a+(1-\theta)b}(\R)).
\]
\end{lemma}

Applying the lemma with successively $(\alpha,\beta,a,b)=(0,1,2,1)$ and $(\alpha,\beta,a,b)=(1,2,1,0)$, we find for any $\theta,\theta'$ in $(0,1)$,
\[
V_2(\R,\R)\subset {\cal C}^{1-\theta}_*(\R,H^{1+\theta}(\R))\cap 
{\cal C}^{2-\theta'}_*(\R,H^{\theta'}(\R)).
\]
Since for any $\epsilon>0$ we have ${\cal C}^{\epsilon}_*\subset{\cal C}^{0}$ and ${\cal C}^{1+\epsilon}_*\subset{\cal C}^{1}$, this concludes the proof of Proposition \ref{th:regularity}.\\

\noindent\textbf{Proof of Lemma \ref{lem:interp}}
It relies on the convexity of the exponential function.
\[
\Vert u\Vert_{{\cal C}^{\theta\alpha+(1-\theta)\beta}_*(\R,H^{\theta a+(1-\theta)b}(\R))}^2
=\sup_{j \ge -1}2^{2j(\theta\alpha+(1-\theta)\beta)}
\sum_{k \ge -1} 2^{2k(\theta a +(1-\theta) b}
\Vert \Delta_j^t \Delta_k^x u\Vert^2 
\]
where $\Delta_j^t$ (resp. $\Delta_j^x$) is the Littlewood Paley operator acting in the time (resp. space) variable.
\[
\begin{array}{lcl}
\ds\sum_{k \ge -1} 2^{2k(\theta a +(1-\theta) b}
\Vert \Delta_j^t \Delta_k^x u\Vert^2 &=&  \ds\sum_{k \ge -1}\Biggl(\Bigl( 2^{2k a}\Vert \Delta_j^t \Delta_k^x u\Vert^2\Bigr)^\theta
\Bigl( 2^{2k b}\Vert \Delta_j^t \Delta_k^x u\Vert^2\Bigr)^{(1-\theta)}\Biggr)\\
&\le &  \Bigl( \ds\sum_{k \ge -1}2^{2k a}\Vert \Delta_j^t \Delta_k^x u\Vert^2\Bigr)^\theta
\Bigl( \ds\sum_{k \ge -1}2^{2k b}\Vert \Delta_j^t \Delta_k^x u\Vert^2\Bigr)^{1-\theta}.
\end{array}
\]
Therefore we have
\[
\Vert u\Vert_{{\cal C}^{\theta\alpha+(1-\theta)\beta}_*(\R,H^{\theta a+(1-\theta)b}(\R))}^2
\le
\Bigl( 
\sup_{j \ge -1}2^{2j\alpha} \ds\sum_{k \ge -1}2^{2k a}\Vert \Delta_j^t \Delta_k^x u\Vert^2
\Bigr)^\theta
\Bigl( 
\sup_{j \ge -1}2^{2j\beta}\sum_{k \ge -1}2^{2k b}\Vert \Delta_j^t \Delta_k^x u\Vert^2
\Bigr)^{1-\theta}
\]
which writes
\[
\Vert u\Vert_{{\cal C}^{\theta\alpha+(1-\theta)\beta}_*(\R,H^{\theta a+(1-\theta)b}(\R))}^2
\le
\Bigl(\Vert u\Vert_{{\cal C}^{\alpha}_*(\R,H^{a}(\R))}^2\Bigr)^\theta\,
\Bigl(\Vert u\Vert_{{\cal C}^{\beta}_*(\R,H^{b}(\R))}^2\Bigr)^{1-\theta}.
\]
\end{proof}
\begin{theorem}\label{th:convswrnl}
  Let $p$ in $H^{2}_0(\Omega)$ and $q$ in $H^{1}_0(\Omega)$. There
  exists a time $T_1 \le T_0^*$ such that for any $T\le T_1$, for any
  initial guess $h_i^\pm$ in $H^1(0,T)$ with the compatibility
  conditions $h_i^+(0)= q(a_{i+1})+p'(a_{i+1})+g^+(p(a_{i+1}))$ and $h_i^-(0)=
  q(a_i)-p'(a_i)+g^-(p(a_{i}))$, the algorithm (\ref{eq:NonoverlappingSchwarz}) is
  defined and converges in $\cup_iV_2(\Omega_i,T)$ to the solution
  $\bar{u}$ of (\ref{clanonlin}).
\end{theorem}
\begin{proof}
We  first prove that the algorithm is well-defined : with the assumptions on $h_i^\pm$ in the theorem, we know by Theorem \ref{th:wpibvpnl} that (\ref{eq:NonoverlappingSchwarz}) defines in each $\Omega_i$  a $u_i^1$
  in $V_2(\Omega_i,T)$, with $\partial_tu_i^1(a_i,\cdot)$, $\partial_tu_i^1(a_{i+1},\cdot)$,
  $\partial_xu_i^1(a_i,\cdot)$ and $\partial_xu_i^1(a_{i+1},\cdot)$ in $H^1(0,T)$ for $T \le T_i^k$. Furthermore, by Lemma \ref{lem:estimfG}, ${\cal B}^-(g^-)u_{i-1}^{1}(a_i,\cdot)$ and  ${\cal B}^+(g^+)u_{i+1}^{1}(a_{i+1},\cdot)$ are in $H^1(0,T)$. As for the compatibility conditions, we have 
\[
{\cal B}^-(g^-)u_{i-1}^{1}(a_i,0)= \lim_{t\rightarrow 0}\bigl(\partial_tu_{i-1}^1(a_i,t)-\partial_xu_{i-1}^1(a_i,t)+g^-(u_{i-1}^{1}(a_i,t))\bigr)
\]
and  by Proposition \ref{th:regularity}, we can pass to the limit and get
\[
{\cal B}^-(g^-)u_{i-1}^{1}(a_i,0)= q(a_i)-p'(a_i)+g^-(p(a_i)).
\]
This, together with the same regularity result on $a_{i+1}$, permits the recursion.\\
 
We define for $T\le \min(T_0^{*},\min_{i}(T^{*k}_i))$, for $k \ge 1$,
the quantities (with $u_i=\bar{u}/_{\Omega_i}$) for $1\le i,j \le I$, $j=i$ or $j=i-1$,
\[
\begin{array}{rll}
e_i^k&=\,\,u_i^k-u_i,\\
f_i^k&=\,\,{\cal F}(u_i^k)-{\cal F}(u_i)\\
% h_{i,j}^{k,-}&= {\cal G}_i^-(u_j^k)-{\cal G}_i^-(u_j),\quad
% h_{i,j}^{k,+}&= {\cal G}_i^+(u_j^k)-{\cal G}_i^+(u_j),
h_{i,j}^{k,-}&= g^-(u_j^k(a_i))- g^-(u_j(a_i)),\quad
h_{i,j}^{k,+}&= g^+(u_j^k(a_{i+1}))-g^+(u_j(a_{i+1})).
\end{array}
\]
The operators ${\cal F}$ is defined in \eqref{eq:SecondMembrek}. The error $e_i^{k}$ in  $\Omega_i$ at iteration $k$ is a solution  of 
\begin{align}
 &  (\partial^2_{t}-\partial^2_x)e_i^{k} + e_i^{k}= f_i^{k}   \mbox{ in $\Omega_i \times (0,T)$}\label{eq:error1}\\
 &  (\partial_{t}-\partial_x)e_i^{k}+ h_{i,i}^{k,-}=(\partial_{t}-\partial_x)e_{i-1}^{k-1} + h_{i,i-1}^{k-1,-}   \mbox{ on }\{a_i\}\times(0,T)\label{eq:error2} \\
 &  (\partial_{t}+\partial_x)e_i^{k}+ h_{i,i}^{k,+}=(\partial_{t}+\partial_x)e_{i+1}^{k-1} + h_{i,i+1}^{k-1,+}  \mbox{ on }\{a_{i+1}\}\times(0,T)\label{eq:error3}  
\end{align}
with vanishing initial values and $e_0^k \equiv 0$, $e_{I+2}^k \equiv 0$, $h_{1,0}^{k,-}\equiv0$, $ h_{I,I+1}^{k,+}=0$ . In order to get a new energy estimate in $\Omega_i$, we multiply (\ref{eq:error1}) by $\partial_te_i^{k}$ and integrate by parts:
\begin{eqnarray}\label{eq:errest1}
\frac{d}{dt}E_{\Omega_i}(e_i^{k})-[\partial_te_i^{k}\partial_xe_i^{k}(a_{i+1},\cdot)-\partial_te_i^{k}\partial_xe_i^{k}(a_{i},\cdot)]=(f_i^{k},\partial_te_i^{k})
\end{eqnarray} 
with $E_{\Omega}(u)=\frac{1}{2}(\|\partial_t u\|^2_{L^2(\Omega)}+\|\partial_x u\|^2_{L^2(\Omega)} + \| u\|^2_{L^2(\Omega)})$. We rewrite the boundary terms using the boundary operators:
\begin{equation}\label{eq:errest2}
\begin{array}{ll}
\partial_te_i^{k}\partial_xe_i^{k}(a_{i+1},\cdot)
=\ds\frac{1}{4}((\partial_{t}+\partial_x)e_i^{k}(a_{i+1},\cdot)+ h_{i,i}^{k,+})^2
\!-\!\frac{1}{4}((\partial_{t}-\partial_x)e_{i}^{k}(a_{i+1},\cdot) + h_{i+1,i}^{k,-})^2 
+R_{i,i+1}^{k}, \\[3mm]
\ds-\partial_te_i^{k}\partial_xe_i^{k}(a_{i},\cdot)
=\ds\frac{1}{4}((\partial_{t}-\partial_x)e_i^{k}(a_{i},\cdot)+ h_{i,i}^{k,-})^2
\!-\!\frac{1}{4}((\partial_{t}+\partial_x)e_{i}^{k}(a_{i},\cdot) + h_{i-1,i}^{k,+})^2
+R_{i,i-1}^{k}.
\end{array}
\end{equation}
The remainders $R_{i,i+1}^{k}$ and $R_{i,i-1}^{k}$ will be evaluated later. We insert (\ref{eq:errest2}) into (\ref{eq:errest1}), and obtain
\[
\begin{split}
\ds\frac{d}{dt}E_{\Omega_i}(e_i^{k})&+ \frac{1}{4}((\partial_{t}-\partial_x)e_{i}^{k}(a_{i+1},\cdot) + h_{i+1,i}^{k,-})^2 + \frac{1}{4}((\partial_{t}+\partial_x)e_{i}^{k}(a_{i},\cdot)+ h_{i-1,i}^{k,+})^2\hspace{4cm}\\[3mm]
&\ds= \frac{1}{4}((\partial_{t}+\partial_x)e_{i}^{k}(a_{i+1},\cdot) + h_{i,i}^{k,+})^2+ 
 \frac{1}{4}((\partial_{t}-\partial_x)e_i^{k}(a_{i},\cdot)+ h_{i,i}^{k,-})^2
+ R_{i,i+1}^{k} +R_{i,i-1}^{k} +(f_i^{k},\partial_te_i^{k}).
\end{split}
\]
Using the transmission conditions (\ref{eq:error2}), (\ref{eq:error3}), we get
\begin{equation}\label{eq:errest3}
\begin{split}
\ds\frac{d}{dt}E_{\Omega_i}(e_i^{k})&+ \frac{1}{4}((\partial_{t}-\partial_x)e_{i}^{k}(a_{i+1},\cdot) + h_{i+1,i}^{k,-})^2 +  \frac{1}{4}((\partial_{t}+\partial_x)e_{i}^{k}(a_{i},\cdot) + h_{i-1,i}^{k,+})^2\hspace{4cm}\\%[3mm]
&= \frac{1}{4}((\partial_{t}+\partial_x)e_{i+1}^{k-1}(a_{i+1},\cdot) + h_{i,i+1}^{k-1,+})^2+ \frac{1}{4}((\partial_{t}-\partial_x)e_{i-1}^{k-1}(a_{i},\cdot) + h_{i,i-1}^{k-1,-})^2\\%[3mm]
&+ R_{i,i+1}^{k} +R_{i,i-1}^{k}+(f_i^{k},\partial_te_i^{k}).
\end{split}
\end{equation}
We sum (\ref{eq:errest3}) on the indexes $i$, $1 \le i \le I$ and integrate in time. We translate the domain indexes in the right-hand side.
Defining 
\[
E_{\partial\Omega_i}(e^k_i)= \ds\frac{1}{4} 
[((\partial_{t}-\partial_x)e_{i}^{k}(a_{i+1},\cdot) + h_{i+1,i}^{k,-})^2
+((\partial_{t}+\partial_x)e_{i}^{k}(a_{i},\cdot) + h_{i-1,i}^{k,+})^2],
\]
we get, since the initial data vanish,
\begin{equation}\label{eq:errest4}
\begin{array}{ll}
\ds\sum_{i=1}^I E_{\Omega_i}(e^{k}_i)(t)&+\ds\int_0^t\sum_{i=1}^I E_{\partial\Omega_i}(e^{k}_i)(s) ds
\le  \ds\int_0^t\sum_{i=1}^I E_{\partial\Omega_i}(e^{k-1}_i)(s) ds\\
& +\ds\int_0^t\sum_{i=1}^{I}\ds( R_{i,i+1}^{k} +R_{i,i-1}^{k})(s)ds
+\ds\int_0^t \sum_{i=1}^I(f_i^{k},\partial_te_i^{k})(s)ds.
\end{array}
\end{equation}
Differentiating the equation and the transmission conditions in time yields the  bound on $\partial_te^k_i$:
\begin{equation}\label{eq:errest5}
\begin{array}{ll}
\ds\sum_{i=1}^I E_{\Omega_i}(\partial_te^{k}_i)(t)&+\ds\int_0^t\sum_{i=1}^I E_{\partial\Omega_i}(\partial_te^{k}_i)(s) ds
\le  \ds\int_0^t\sum_{i=1}^I E_{\partial\Omega_i}(\partial_te^{k-1}_i)(s) ds\\
& +\ds\int_0^t\sum_{i=1}^{I}\ds(  \tilde{R}_{i,i+1}^{k} + \tilde{R}_{i,i-1}^{k})(s)ds
+\ds\int_0^t \sum_{i=1}^I(\partial_tf_i^{k+1},\partial_{tt}e_i^{k})(s)ds.
\end{array}
\end{equation}
We now estimate the remainders. We start with $R_{i,i+1}^{k}$ (ignoring the superscript $k$):
\[
R_{i,i+1}=  \frac{1}{4}[- h_{i,i}^+(2(\partial_{t}e_{i}+\partial_xe_{i})(a_{i+1},\cdot)+h_{i,i}^+)
         + h_{i+1,i}^-(2(\partial_{t}e_{i}-\partial_xe_{i})(a_{i+1},\cdot)+h_{i+1,i}^-) ],
\]
and we get a bound on the integral of $ R_{i,i+1}$:
\[
\begin{array}{ll}
\ds\int_0^t R_{i,i+1}(s) ds\le  \frac{3}{4}( \|h_{i,i}^+\|^2_{H^{\frac{1}{2}}(0,t)}+\|h_{i+1,i}^-\|^2_{H^{\frac{1}{2}}(0,t)})+ \frac{1}{2}(\|\partial_{t}e_{i}(a_{i+1},\cdot)\|^2_{H^{-\frac{1}{2}}(0,t)}+\|\partial_{x}e_{i}(a_{i+1},\cdot)\|^2_{H^{-\frac{1}{2}}(0,t)}).
\end{array}
\]
We can treat $R_{i,i-1}$,$\tilde{R}_{i,i-1}$, and $\tilde{R}_{i,i-1}$ the same way and obtain
\[
\begin{array}{ll}
\ds\int_0^t R_{i,i-1}(s) ds\le  & \frac{3}{4}( \|h_{i,i}^-\|^2_{H^{\frac{1}{2}}(0,t)}+\|h_{i-1,i}^+\|^2_{H^{\frac{1}{2}}(0,t)})
+ \frac{1}{2}(\|\partial_{t}e_{i}(a_{i},\cdot)\|^2_{H^{-\frac{1}{2}}(0,t)}+\|\partial_{x}e_{i}(a_{i},\cdot)\|^2_{H^{-\frac{1}{2}}(0,t)}),
\end{array}
\]
\[
\begin{array}{ll}
\ds\int_0^t \tilde{R}_{i,i+1}(s) ds\le  & \frac{3}{4}( \|\partial_th_{i,i}^+\|^2_{H^{\frac{1}{2}}(0,t)}+\|\partial_th_{i+1,i}^-\|^2_{H^{\frac{1}{2}}(0,t)})
+ \frac{1}{2}(\|\partial_{t}^2e_{i}(a_{i+1},\cdot)\|^2_{H^{-\frac{1}{2}}(0,t)}+\|\partial_{xt}e_{i}(a_{i+1},\cdot)\|^2_{H^{-\frac{1}{2}}(0,t)}),
\end{array}
\]
\[
\begin{array}{ll}
\ds\int_0^t \tilde{R}_{i,i-1}(s) ds\le  &  \frac{3}{4}( \|\partial_th_{i,i}^-\|^2_{H^{\frac{1}{2}}(0,t)}+\|\partial_th_{i-1,i}^+\|^2_{H^{\frac{1}{2}}(0,t)})
+ \frac{1}{2}(\|\partial_{t}^2e_{i}(a_{i},\cdot)\|^2_{H^{-\frac{1}{2}}(0,t)}+\|\partial_{xt}e_{i}(a_{i},\cdot)\|^2_{H^{-\frac{1}{2}}(0,t)}).
\end{array}
\]
At point $a_i$ for instance, by the Trace Theorem, there is a  constant $C_3$ independent of $T$, such that for any $\alpha$ with $|\alpha|=1$, we have
\[
 \begin{array}{c}
\|\partial^\alpha e_{i}(a_{i},\cdot)\|_{H^{-\frac{1}{2}}(0,t)} \le  \|\partial^\alpha e_{i}(a_{i},\cdot)\|_{L^{2}(0,t)} \le C_3  \|\partial^\alpha e_{i}\|_{H^1(\Omega_i\times(0,t))},\\[3mm] 
\|\partial^\alpha \partial_{t}e_{i}(a_{i},\cdot)\|_{H^{-\frac{1}{2}}(0,t)} \le  \|\partial^\alpha e_{i}(a_{i},\cdot)\|_{H^{\frac{1}{2}}(0,t)} \le C_3  \|\partial^\alpha e_{i}\|_{H^1(\Omega_i\times(0,t))}, 
\end{array}
\]
which gives our first bounds on the remainders:
\[
\ds\int_0^t R_{i,i+1}(s) ds\le \frac{1}{2}( \|h_{i,i}^+\|^2_{H^{\frac{1}{2}}(0,t)}+\|h_{i+1,i}^-\|^2_{H^{\frac{1}{2}}(0,t)})+ \frac{C_3^2}{2}\|e_i\|^2_{H^2(\Omega_i\times(0,t))},
\]
\[
\ds\int_0^t \tilde{R}_{i,i+1}(s) ds\le \frac{1}{2}( \|\partial_th_{i,i}^+\|^2_{H^{\frac{1}{2}}(0,t)}+\|\partial_th_{i+1,i}^-\|^2_{H^{\frac{1}{2}}(0,t)})+ \frac{C_3^2}{2}\|e_i\|^2_{H^2(\Omega_i\times(0,t))}.
\]
We now   insert the previous estimates in  (\ref{eq:errest4}) and (\ref{eq:errest5}). By Cauchy-Schwarz inequality we get
\begin{multline}\label{eq:errest6}
\ds\sum_{i=1}^I E_{\Omega_i}(e^{k}_i)(t)+\ds\int_0^t\sum_{i=1}^I E_{\partial\Omega_i}(e^{k}_i)(s) ds
\le \\ 
\ds\int_0^t\sum_{i=1}^I E_{\partial\Omega_i}(e^{k-1}_i)(s) ds
+ \frac{C_3^2+1}{2}\sum_{i=1}^I\|e_i^{k}\|^2_{H^2(\Omega_i\times(0,t))}
+\ds\frac{1}{2}\sum_{i=1}^I\int_0^t\|f_i^{k}\|_{L^2(\Omega_i)}^2(s)ds \\
+ \ds\frac{3}{4}\sum_{i=1}^{I}\ds( \|h_{i,i}^{k,+}\|^2_{H^{\frac{1}{2}}(0,t)}+\|h_{i+1,i}^{k,-}\|^2_{H^{\frac{1}{2}}(0,t)})
\ds+\frac{3}{4}\ds\sum_{i=1}^{I}\ds( \|h_{i,i}^{k,-}\|^2_{H^{\frac{1}{2}}(0,t)}+\|h_{i-1,i}^{k,+}\|^2_{H^{\frac{1}{2}}(0,t)})
 ,
\end{multline}
\begin{multline}\label{eq:errest7}
\ds\sum_{i=1}^I E_{\Omega_i}(\partial_te^{k}_i)(t)+\ds\int_0^t\sum_{i=1}^I E_{\partial\Omega_i}(\partial_te^{k}_i)(s) ds
\le\\
\ds\int_0^t\sum_{i=1}^I E_{\partial\Omega_i}(\partial_te^{k-1}_i)(s) ds
+ \frac{C_3^2+1}{2}\sum_{i=1}^I\|e_i^{k}\|^2_{H^2(\Omega_i\times(0,t))}
+\frac{1}{2}\ds \sum_{i=1}^I\int_0^t\|\partial_tf_i^{k}\|_{L^2(\Omega_i)}^2 (s)ds 
\\
+ \ds\frac{3}{4}\sum_{i=1}^{I}\ds( \|\partial_th_{i,i}^{k,+}\|^2_{H^{\frac{1}{2}}(0,t)}+\|\partial_th_{i+1,i}^{k,-}\|^2_{H^{\frac{1}{2}}(0,t)})
+\frac{3}{4}\ds\sum_{i=1}^{I}\ds( \|\partial_th_{i,i}^{k,-}\|^2_{H^{\frac{1}{2}}(0,t)}+\|\partial_th_{i-1,i}^{k,+}\|^2_{H^{\frac{1}{2}}(0,t)})\\
\end{multline}
Adding (\ref{eq:errest6}) and (\ref{eq:errest7}), we can write
\[
\begin{array}{l}
\ds\sum_{i=1}^I (E_{\Omega_i}(e^{k}_i) +E_{\Omega_i}(\partial_te^{k}_i))(t)+\ds\int_0^t\sum_{i=1}^I (E_{\partial\Omega_i}(e^{k}_i)+E_{\partial\Omega_i}(\partial_te^{k}_i))(s) ds\\
\le  \ds\int_0^t\sum_{i=1}^I (E_{\partial\Omega_i}(e^{k-1}_i)+ E_{\partial\Omega_i}(\partial_te^{k-1}_i))(s) ds + (C_3^2+1)\ds\sum_{i=1}^I\|e_i^{k}\|^2_{H^2(\Omega_i\times(0,t))}\\
+\ds\frac{3}{4}\sum_{i=1}^{I}\ds( \|h_{i,i}^{k ,+}\|^2_{H^{\frac{1}{2}}(0,t)}+\|h_{i+1,i}^{k ,+}\|^2_{H^{\frac{1}{2}}(0,t)}+ \|h_{i,i}^{k ,-}\|^2_{H^{\frac{1}{2}}(0,t)}+\|h_{i-1,i}^{k ,+}\|^2_{H^{\frac{1}{2}}(0,t)})\\
 + \ds\frac{3}{4}\sum_{i=1}^{I}\ds( \|\partial_th_{i,i}^{k ,+}\|^2_{H^{\frac{1}{2}}(0,t)}+\|\partial_th_{i+1,i}^{k ,+}\|^2_{H^{\frac{1}{2}}(0,t)}+ \|\partial_th_{i,i}^{k ,-}\|^2_{H^{\frac{1}{2}}(0,t)}+\|\partial_th_{i-1,i}^{k ,+}\|^2_{H^{\frac{1}{2}}(0,t)})\\
+\ds \sum_{i=1}^I\frac{1}{2}(\|f_i^{k }\|^2_{L^2(\Omega\times(0,T))}+\|\partial_tf_i^{k }\|^2_{L^2(\Omega\times(0,T))}).
\end{array}
\]
We now estimate the quantities involving the $h_{i,j}^{k,\pm}$. Refining the results in Lemma \ref{lem:estimfG}, we have a real positive increasing function $\theta_2$,  such that
\[
\ds \|h_{i,j}^{k,\pm}\|^2_{H^{\frac{1}{2}}(0,t)}+ \|\partial_th_{i,j}^{k,\pm}\|^2_{H^{\frac{1}{2}}(0,t)}\le \theta_2^2(\sum_{\stackrel{|\alpha|\le 2}{\alpha \ne (0,2)}
}(\|\partial^\alpha u_i\|_{L^2(\Omega\times(0,T))}^2+\|\partial^\alpha u_i^k\|_{L^2(\Omega\times(0,T))}^2))\sum_{\stackrel{|\alpha|\le 2}{\alpha \ne (0,2)}
}\|\partial^\alpha e_i^k\|_{L^2(\Omega\times(0,T))}^2
\]
which gives
\begin{equation}\label{eq:errest8}
\begin{array}{l}
 \ds\sum_{i=1}^I (E_{\Omega_i}(e^{k}_i) +E_{\Omega_i}(\partial_te^{k}_i))(t)+\ds\int_0^t\sum_{i=1}^I (E_{\partial\Omega_i}(e^{k}_i)+E_{\partial\Omega_i}(\partial_te^{k}_i))(s) ds\\
\hspace{2cm}\le  \ds\int_0^t\sum_{i=1}^I (E_{\partial\Omega_i}(e^{k-1}_i)+ E_{\partial\Omega_i}(\partial_te^{k-1}_i))(s) ds \\
\hspace{2cm}+\ds \sum_{i=1}^I\theta_3(\sum_{\stackrel{|\alpha|\le 2}{\alpha \ne (0,2)}
}(\|\partial^\alpha u_i\|_{L^2(\Omega\times(0,T))}^2+\|\partial^\alpha u_i^k\|_{L^2(\Omega\times(0,T))}^2))\ds\|e_i^{k}\|^2_{H^2(\Omega_i\times(0,t))}\\
\hspace{2cm}+\ds \sum_{i=1}^I\frac{1}{2}(\|f_i^{k}\|^2_{L^2(\Omega\times(0,T))}+\|\partial_tf_i^{k}\|^2_{L^2(\Omega\times(0,T))}),
\end{array}
\end{equation}
with $\theta_3= 3\theta_2^2+ C_3^2+1$.  We now evaluate the  terms in the right-hand side. We first note that 
\[
\|e^{k}_i\|^2_{H^2(\Omega_i\times(0,t))}\le \int_0^t(E_{\Omega_i}(e^{k}_i)+E_{\Omega_i}(\partial_te^{k}_i)+\|\partial_{xx}e^{k}_i\|_{L^2(\Omega_i)}^2)ds,
\] 
and evaluate $\|\partial_{xx}e^{k}_i\|^2_{L^2(\Omega_i\times(0,t))}$ by  equation (\ref{eq:error1}):
\[
\|\partial_{xx}e^{k}_i\|^2_{L^2(\Omega_i\times(0,t))}\le 3(\|\partial_{tt}e^{k}_i\|^2_{L^2(\Omega_i\times(0,t))} +\|e^{k}_i\|^2_{L^2(\Omega_i\times(0,t))}+ \|f^{k}_i\|^2_{L^2(\Omega_i\times(0,t))}),
\] 
from which we deduce
\[
\|e^{k}_i\|^2_{H^2(\Omega_i\times(0,t))}\le 4\int_0^t(E_{\Omega_i}(e^{k}_i)+E_{\Omega_i}(\partial_te^{k}_i))ds + 3 \|f^{k}_i\|^2_{L^2(\Omega_i\times(0,t))}.
\] 
There remains only in (\ref{eq:errest8})
\begin{equation}\label{eq:errest9}
\begin{array}{l}
 \ds\sum_{i=1}^I (E_{\Omega_i}(e^{k}_i) +E_{\Omega_i}(\partial_te^{k}_i))(t)+\ds\int_0^t\sum_{i=1}^I (E_{\partial\Omega_i}(e^{k}_i)+E_{\partial\Omega_i}(\partial_te^{k}_i))(s) ds\\
\hspace{1cm}\le  \ds\int_0^t\sum_{i=1}^I (E_{\partial\Omega_i}(e^{k-1}_i)+ E_{\partial\Omega_i}(\partial_te^{k-1}_i))(s) ds \\
\hspace{1cm}+\ds \sum_{i=1}^I 4\theta_3(\sum_{\stackrel{|\alpha|\le 2}{\alpha \ne (0,2)}
}(\|\partial^\alpha u_i\|_{L^2(\Omega\times(0,T))}^2+\|\partial^\alpha u_i^k\|_{L^2(\Omega\times(0,T))}^2))\ds\int_0^t(E_{\Omega_i}(e^{k}_i)+E_{\Omega_i}(\partial_te^{k}_i))ds \\
+\ds \sum_{i=1}^I (3 \theta_3(\sum_{\stackrel{|\alpha|\le 2}{\alpha \ne (0,2)}
}(\|\partial^\alpha u_i\|_{L^2(\Omega\times(0,T))}^2+\|\partial^\alpha u_i^k\|_{L^2(\Omega\times(0,T))}^2))+ \frac{1}{2})(\|f_i^{k}\|^2_{L^2(\Omega\times(0,T))}+\|\partial_tf_i^{k}\|^2_{L^2(\Omega\times(0,T))}).
\end{array}
\end{equation}
Again, as in Lemma \ref{lem:estimfG}, there exists a positive increasing function $\theta_4$ such that 
\[
\ds\|f_i^k(t,\cdot)\|_{L^2(\Omega)}^2+ \|\partial_tf_i^k(t,\cdot)\|_{L^2(\Omega)}^2\le  \theta_4^2( \sum_{\stackrel{|\alpha|\le 2}{\alpha \ne (0,2)}
}(\|\partial^\alpha u_i\|_{L^2(\Omega)}^2+\|\partial^\alpha u_i^k\|_{L^2(\Omega)}^2)\sum_{\stackrel{|\alpha|\le 2}{\alpha \ne (0,2)}
}\|\partial^\alpha e_i^k\|_{L^2(\Omega)}^2,
\] 
the latter sum means that no term $\partial_{xx}$ are present, therefore we can bound the sum by twice the energy. Furthermore we know that the energy of $\bar{u}$ is bounded on the interval $(0,T)$. Thus  there exists a new positive increasing function $\theta_5$, depending on $u$, such that 
\begin{equation}\label{eq:boundf}
\ds\|f_i^k(t,\cdot)\|_{L^2(\Omega)}^2+ \|\partial_tf_i^k(t,\cdot)\|_{L^2(\Omega)}^2\le (\theta_5^2(
E_{\Omega_i}(e_i^k) +E_{\Omega_i}(\partial_te^{k}_i))(E_{\Omega_i}(e_i^k)+E_{\Omega_i}(\partial_te^{k}_i)))(t).
\end{equation}
We insert (\ref{eq:boundf}) into (\ref{eq:errest9}), and get (with a new function $\theta_6$)
\begin{equation}\label{eq:errest10}
\begin{array}{l}
 \ds\sum_{i=1}^I (E_{\Omega_i}(e^{k}_i) +E_{\Omega_i}(\partial_te^{k}_i))(t)+\ds\int_0^t\sum_{i=1}^I (E_{\partial\Omega_i}(e^{k}_i)+E_{\partial\Omega_i}(\partial_te^{k}_i))(s) ds\\
\hspace{1cm}\le  \ds\int_0^t\sum_{i=1}^I (E_{\partial\Omega_i}(e^{k-1}_i)+ E_{\partial\Omega_i}(\partial_te^{k-1}_i))(s) ds \\
\hspace{1cm}+\ds  \sum_{i=1}^I \int_0^t (\theta_6(
E_{\Omega_i}(e_i^k) +E_{\Omega_i}(\partial_te^{k}_i)))(E_{\Omega_i}(e_i^k)+E_{\Omega_i}(\partial_te^{k}_i))(s)ds.
\end{array}
\end{equation}
For clarity we define 
\[
\begin{array}{l}
E_{int}^k=\ds\sum_{i=1}^I (E_{\Omega_i}(e^{k}_i) +E_{\Omega_i}(\partial_te^{k}_i)),\quad
E_{b}^k=\ds\sum_{i=1}^I (E_{\partial\Omega_i}(e^{k}_i) +E_{\partial\Omega_i}(\partial_te^{k}_i)).\\
\end{array}
\] 
and we can rewrite (\ref{eq:errest10}) as
\begin{equation}\label{eq:errest11}
\ds E_{int}^{k}(t)+ \int_0^t E_{b}^{k}(s)ds \le \int_0^t E_{b}^{k-1}(s)ds+  \int_0^t\theta_6(E_{int}^{k}(s)) E_{int}^{k}(s)ds
\end{equation}
Summing in $k$, we define $\widetilde{E}_{int}^K=\sum_{k=1}^K  E_{int}^{k}$, and we have 
\[
\ds \widetilde{E}_{int}^{K}(t)+  \int_0^t E_{b}^{K}(s)ds \le   \int_0^t E_{b}^{0}(s)ds+  \int_0^t\theta_6(\widetilde{E}_{int}^{K}(s)) \widetilde{E}_{int}^{K}(s)ds.
\] 
Let now $C>0$. If $t$ tends to $0$, $ \int_0^t E_{b}^{1}(s)ds+ tC\theta_6(C)$ tends to $0$. Therefore there exists a  $T_1$  such that    $ \int_0^{T_1} E_{b}^{1}(s)ds+ T_1C\theta_6(C)=C$, and so we have for $t \le T_1$, 
\[
\ds \widetilde{E}_{int}^{K}(t)\le   C.
\]
We conclude that $u_i^k$ exists on the time interval $(0,T_1)$, and that $\sum_{i=1}^I(E_{\Omega_i}(e^{k}_i) +E_{\Omega_i}(\partial_te^{k}_i))$ tends to $0$ when $k$ tends to infinity: the sequence $u_i^k$ converges to $\bar{u}$ on $(0,T_1)$ in each subdomain in the norm of energy.
\end{proof}

%=================================================================
\section{A Finite Volume Discretization}
\label{sec:fv}
%=================================================================
We use here a finite volumes scheme, which has been described in
\cite{GHN:2003:OSW} for the linear one-dimensional wave equation, and
extended to the non linear boundary value problems in the frame of
absorbing boundary conditions in
\cite{Sz:1900:ABC}. We restrict ourselves to uniform meshes in time and
space.
%--------------------------------------------------------------
\subsection{Discretization of the Subdomain Problem (\ref{eq:IBVPNL})}
\label{subsec:disc}
%--------------------------------------------------------------
%
The domain $\Omega\times (0,T)$ is meshed by a rectangular grid, with
uniform mesh sizes $\Delta x$ and $\Delta t$.  There are $J+1$ points
in space with $\Delta x=(a_{+}-a_-)/J$, and $N+1$ points in time, with
$\Delta t=T/N$. We denote the numerical approximation to
$u(a_-+j\Delta x,n\Delta t)$ by $U(j,n)$.  We introduce the notations:
{\small%
\begin{equation}\label{eq:finiteDifferences}
  \arraycolsep0.2em
  \begin{array}{rlclc}
    \ds \dtp U(j,n) & =  \ds\frac{U(j,n+1)-U(j,n)}{\Delta t},&
    \ds \dtm U(j,n) & =  \ds\frac{U(j,n)-U(j,n-1)}{\Delta t},\\[3mm]
    \ds \dxp U(j,n) & =  \ds\frac{U(j+1,n)-U(j,n)}{\Delta x},&
    \ds \dxm U(j,n) & =  \ds\frac{U(j,n)-U(j-1,n)}{\Delta x},\\[3mm]
    \ds \dxz U(j,n) & =  \ds\frac{U(j+1,n)-U(j-1,n)}{2\Delta x},&
    \ds \dtz U(j,n) & =  \ds\frac{U(j,n+1)-U(j,n-1)}{2\Delta t}, \\[3mm]
    \ds \dtmm U(j,n) & =  \ds\frac{3\,U(j,n)-4U(j,n-1)+U(j,n-2)}{2\Delta t},&
    \ds \dtms U(j,n) & =  \left\{
                           \begin{array}{l}
                           \dtmm U(j,n) \mbox{ for }n \ge 2,\\
                            \dtm U(j,n)  \mbox{ for }n =1.
                            \end{array}
                            \right.
  \end{array}
\end{equation}
} 
The last finite derivative in (\ref{eq:finiteDifferences}) is a second
order approximation of $\partial_t u$, to be used in the nonlinear
term, in order to design an explicit scheme.

The scheme in the interior writes
\begin{equation}
  \label{eq:interiornlscheme}
  \left(\dtp\dtm-\dxp\dxm\right)U(j,n)-f(U(j,n),\dtms U(j,n),\dxz U(j,n))=0 \ \mbox{ in } \ie{1,J-1}\times\ie{1,N-1}. 
\end{equation}
We define the discrete initial value as
\begin{equation*}\label{eq:discreteinitcond}
P(j)=p(a_-+j\Delta x), \quad Q(j)=q(a_-+j\Delta x),
\end{equation*}
and  we obtain the initial scheme
\begin{equation}\label{eq:schemeinit}
  (\dtp - \frac{\Delta t}{2}\dxp\dxm)U(j,0)= Q(j)+ \frac{\Delta t}{2}f(P(j),Q(j),\dxz Q(j)),\quad \mbox{ in } \ie{1,J-1}.
\end{equation}
For the boundary conditions, we define the discrete boundary operators
as:
{\footnotesize
\begin{equation}  \label{eq:bm}
  \begin{array}{l}
 \ds B^{-}(f,g^-)\,U(0,n)=\ds (\dtz-\dxp+\frac{\Delta x}{2}\dtp\dtm) U(0,n)
\ds- \frac{\Delta x}{2}f(U(0,n),\dtms U(0,n),\dxp U(0,n)) + g^-(U(0,n)),\\[2mm]
  B^-(f,g^-)\,U(0,0)=\ds(\dtp-\dxp+\frac{\Delta x}{\Delta t}\dtp)\,U(0,0)-\frac{\Delta x}{\Delta t}Q(0)
-\frac{\Delta x}{2}f(P(0),Q_i(0),\dxp P(0))+g^-(P(0)), 
\end{array}
\end{equation}
}
{\footnotesize
\begin{equation}  \label{eq:bp}
  \begin{array}{l}
\ds B^{+}(f,g^+)\,U(J,n)=\ds (\dtz+\dxm+\frac{\Delta x}{2}\dtp\dtm) U(J,n)
- \frac{\Delta x}{2}f(U(J,n),\dtms U_i(J,n),\dxm U(J,n)) + g^+(U(J,n)).\\[2mm]
B^+(f,g^+)\,U(J,0)=\ds(\dtp+\dxm+\frac{\Delta x}{\Delta t}\dtp)U(J,0)-\frac{\Delta x}{\Delta t}Q(J)
-\frac{\Delta x}{2}f(P(J),Q_i(J),\dxm P(J)) +\ds g^+(P(J)).
\end{array}
\end{equation}
}%
We define the boundary data as 
\begin{equation}\label{eq:discretebdata}
H^\pm(n)=\frac{1}{\Delta t}\int_{t_n-\Delta t/2}^{t_n+\Delta t/2} h^\pm(\tau)d\tau, \quad 1 \le n \le N,H^\pm(0)= \frac{2}{\Delta t}\int_{0}^{\Delta t/2} h^\pm(\tau)d\tau
\end{equation}
The discretization of Problem \eqref{eq:IBVPNL} is now  given by \eqref{eq:interiornlscheme}, \eqref{eq:schemeinit}, with boundary conditions 
 \begin{equation}\label{eq:NumericalSchemebc}
   B^-(f,g^-)U (0,\cdot)  =H^-\ ,\   B^+(f,g^+)U(J,\cdot) = H^+\mbox{ in } \,  \ie{0,n \le N},
\end{equation}
where the discrete boundary operators $B^\pm$ are given in
(\ref{eq:bm}), (\ref{eq:bp}). \\

Our numerical computations indicate that this scheme is second order both in space and time.

%--------------------------------------------------------------
\subsection{The Discrete Schwarz Waveform Relaxation Algorithm}
\label{subsec:dswr}
%--------------------------------------------------------------
The equation
is now discretized on each subdomain $\Omega_i \times (0,T)$,
$i=1,\ldots,I$ separately, using an  uniform mesh with sizes $\Delta x$ and $\Delta t$. There are $J_i+1$ points in space and $N+1$ grid points in time in subdomain
$\Omega_i$, with $\Delta
x=(a_{i+1}-a_i)/J_i$ and $\Delta t=T/N$. We denote the numerical approximation
to $u_i^k(a_i+j\Delta x,n\Delta t)$ on $\Omega_i$ at iteration step
$k$ by $U_i^k(j,n)$. 

The problem in domain $\Omega_i$ is now defined through boundary data $H^{\pm,\, k-1}_i$, coming from the neighboring subdomains $\Omega_{i\pm 1}$ .Therefore, we define the extraction operator from domain $\Omega_{i}$ to his neighbours as:
{\footnotesize
\begin{equation}  \label{eq:tildebp}
\begin{array}{l}
  \widetilde{B}^{+}(f,g^+)\,U_{i}(0,n)=\ds(\dtz +\dxp -\frac{\Delta x}{2}\dtp\dtm )U_{i}(0,n) 
  +\frac{\Delta x}{2}f(U_{i}(0,n),\dtms U_{i}(0,n),\dxp U_{i}(0,n))\ds+ g^+(U_{i}(0,n)),\\[2mm]
  \widetilde{B}^+(f,g^+)\,U_{i}(0,0)=\ds(\dtp+\dxp-\frac{\Delta x}{\Delta t}\dtp)\,U_{i}(0,0)+\frac{\Delta x}{\Delta t}Q_{i}(0)
  +\frac{\Delta x}{2}f(P_{i}(0),Q_{i}(0),\dxp P_{i}(0))+g^+(P_{i}(0)).
   \end{array}
\end{equation}
}
{\footnotesize
\begin{equation}  \label{eq:tildebm}
  \begin{array}{l}
\ds \widetilde{B}^{-}(f,g^-)U_{i}(\ji,n)= \ds (\dtz-\dxm-\frac{\Delta x}{2}\dtp\dtm) U_{i}(\ji,n)
+ \frac{\Delta x}{2}f(U_{i}(\ji,n),\dtms U_{i}(\ji,n),\dxm U_{i}(\ji,n))\ds + g^-(U_{i}(\ji,n)),\\[2mm]
\widetilde{B}^-(f,g^-)\,U_{i}(\ji,0)=\ds(\dtp-\dxm-\frac{\Delta x}{\Delta t}\dtp)\,U_{i}(\ji,0)+\frac{\Delta x}{\Delta t}Q_{i}(\ji)
+\frac{\Delta x}{2}f(P_{i}(\ji),Q_{i}(\ji),\dxm P_{i}(\ji))+g^-(P_{i}(\ji)). 
\end{array}
\end{equation}
}
The discrete Schwarz waveform relaxation algorithm on
subdomains $\Omega_i$, $i=1,\ldots, I$ is defined as follows. An initial guess  $\{H^{0,\pm}_i\}_{1\le i \le I}$ is given. For $k \ge 1$, we solve  
\begin{equation}\label{eq:NumericalSchemeSWR}
  \begin{array}{l} 
    \left[
      \begin{array}{l} 
        \left(\dtp\dtm-\dxp\dxm\right)U_i^k=f(U_i^k,\dtms U_i^k,\dxz U_i^k) 
       \mbox{ in } \,   \ie{1,J_i-1} \times  \ie{1,N},\\[3mm]
       \ds   U_i^k(\cdot,0)=P_i,\quad (\dtp - \frac{\Delta t}{2}\dxp\dxm)U^k_i(\cdot,0)= Q_i+ \frac{\Delta t}{2}f(P_i,Q_i,\dxz Q_i)\mbox{ in } \,  \ie{1,J_i-1},\\[3mm]
       B^-(f,g^-)U^{k}_i (0,\cdot) =  H^{-,\, k-1}_{i}\ ,\ 
       B^+(f,g^+)U^{k}_i (J_i,\cdot) =  H^{+,\, k-1}_{i} \mbox{ in } \, \ie{0,n \le N}\\[3mm]
     \end{array}
   \right.
   \\[5mm]
\ds H_{i}^{-,\, k}= \widetilde{B}^-(f,g^-)U^{k}_{i-1}(J_{i-1},\cdot)\ ,\ H_{i}^{+,\, k} = \widetilde{B}^+(f,g^+)U^{k}_{i+1}(0,\cdot) \mbox{ in } \, \ie{0,n \le N}.
\end{array}
\end{equation}
As in the continuous algorithm, we set $U^k_0 \equiv 0$ and $U^k_{I+2}
\equiv 0$.  

We denote by $\bar{U}$ the discrete approximation of problem
(\ref{clanonlin}), obtained by solving
(\ref{eq:interiornlscheme},\ref{eq:schemeinit},\ref{eq:NumericalSchemebc})
on $\Omega=(a,b)$ with $J=\sum_{i=1}^I J_i$ intervals of length
$\Delta x$, and $H^\pm=0$.  Each subproblem is an explicit scheme, thus has a unique
solution.  Therefore the Schwarz waveform Relaxation Algorithm is well
defined. If it converges, the limit in each subdomain is denoted by
$V_i$. It satisfies the same scheme as $\bar{U}$ at initial time, in the
interior and on the exterior boundaries. At point $a_i$ it satisfies for any $n \ge 0$
\begin{equation}\label{eq:tclim}
B^-(f,g^-)V_i (0,n)= \widetilde{B}^-(f,g^-)V_{i-1}(J_{i-1},n)\ ,\  B^+(f,g^+)V_{i-1} (J_{i-1},n) = \widetilde{B}^+(f,g^+)V_{i}(0,n).
\end{equation}
\begin{theorem}\label{th:conv}
Suppose that  $f$ is affine in the third variable $\partial_x u$. Then, if the discrete algorithm converges, it converges to the discrete approximation $\bar{U}$ of problem
(\ref{clanonlin}).
\end{theorem}
\begin{proof} 
Since $p$ is in $H^2(\Omega)$ and $q$  is in $H^1(\Omega)$, there are both continuous and we have for any $i$ $V_i(0,0)=P_i(0)=V_{i-1}(\jim,0)=P_{i-1}(\jim)=p(a_i)$, and  $Q_{i-1}(\jim)=Q_i(0)=q(a_i)$. We  write the transmission conditions (\ref{eq:tclim}) for $n=0$. The nonlinear terms containing $g^\pm$ on both sides cancel out, and we have 
\begin{multline}\label{eq:lim1}
(\dtp-\dxp+ \frac{\Delta x}{\Delta t}\dtp)V_i(0,0)
-\frac{\Delta x}{\Delta t}q(a_i) - \frac{\Delta x}{2}f(p(a_i),q(a_i),\dxp p(a_i) )
=\\
(\dtp-\dxm- \frac{\Delta x}{\Delta t}\dtp)V_{i-1}(J_{i-1},0)
+\frac{\Delta x}{\Delta t}q(a_i) + \frac{\Delta x}{2}f(p(a_i),q(a_i),\dxm p(a_i) ),
\end{multline}
\begin{multline}\label{eq:lim2}
(\dtp+\dxp-\frac{\Delta x}{\Delta t}\dtp)V_i(0,0)
+\frac{\Delta x}{\Delta t}q(a_i) + \frac{\Delta x}{2}f(p(a_i),q(a_i),\dxp p(a_i) )
=\\
(\dtp+\dxm+ \frac{\Delta x}{\Delta t}\dtp)V_{i-1}(J_{i-1},0)
-\frac{\Delta x}{\Delta t}q(a_i) - \frac{\Delta x}{2}f(p(a_i),q(a_i)),\dxm p(a_i)).
\end{multline}
Adding \eqref{eq:lim1} and  \eqref{eq:lim2} yields $\dtp V_i(0,0)=\dtp V_{i-1}(J_{i-1},0)$, and hence  $V_i(0,1)=V_{i-1}(J_{i-1},1)$. We define now $\tilde{V}(j,n)$ for $0 \le j \le J$ as $\tilde{V}(\tilde{j},n)=V_i(j,n)$ if $\tilde{j}\Delta x= a_i + j\Delta x$, with $1 \le j \le J_i$. With the assumption on $f$, since $\dxz=(\dxp+\dxm)/2$,  we can rewrite \eqref{eq:lim1} as
\[
(\dxp-\dxm -2\frac{\Delta x}{\Delta t}\dtp)\tilde{V}(\tilde{j}_i,0)
+2\frac{\Delta x}{\Delta t}q(a_i) + \Delta x f(p(a_i),q(a_i),\dxz p(a_i))
=0,
\]
with $\tilde{j}_i=J_1+\dotsc +J_i$, which, multiplying by $-\frac{\Delta t }{2\Delta x}$,  proves that $ \tilde{V}$ is solution of (\ref{eq:schemeinit}) at any point, and therefore $U$ and $ \tilde{V}$ coincide at time 0 and 1. A simple recursion with the explicit schemes now proves that,  for any $n$, for any $i$, $V_i(0,n)=V_{i-1}(\jim,n)$, and therefore  $U=\tilde{V}$.  
\end{proof}

\begin{remark}
The assumption on $f$ in Theorem \ref{th:conv} is fulfilled when $f(u,u_t,u_x)=f_1(u)+f_2(u)u_t+f_3(u)u_x$. 
\end{remark}
%=================================================================
\section{Convergence of the Discrete Algorithm}
\label{sec:cd}
%=================================================================
According to Theorem \ref{th:conv}, we suppose that $f$ is affine in $\partial_x u$. We introduce the linear transmission operators defined by
\begin{equation*}
  T^\pm=B^\pm(0,0),\quad \ds\widetilde{T}^\pm=\widetilde{B}^\pm(0,0)
\end{equation*} 
We note $\bar{U}_i=\bar{U}/_{\Omega_i}$. With these notations, the error $\bar{U}_i^k=U_i^k-\bar{U}_i$ is solution of the linear problem
\begin{align}
  \left(\dtp\dtm-\dxp\dxm+1\right)\bar{U}_i^k=F_i^k,  \quad  \mbox{ in } \,   \ie{1,J_i-1} \times  \ie{1,N}\label{eq:errorSWRint}\\[3mm]
  \end{align}
with the initial value $\bar{U}_i^k(j,1)=\bar{U}_i^k(j,0)=0$ and the transmission conditions
 \begin{equation}\label{eq:errorSWRbc}
\begin{aligned}
                T^-\bar{U}^{k}_i (0,\cdot) + G_i^{-,\, k}= 
           \widetilde{T}^-\bar{U}^{k-1}_{i-1}(J_{i-1},\cdot)+ \widetilde{G}_{i-1}^{-,\, k-1}, \quad \mbox{ in } \,  \ie{0,N},\\ 
              T^+\bar{U}^{k}_i (J_i,\cdot) + G_i^{+,\, k}= 
      \widetilde{T}^+\bar{U}^{k-1}_{i+1}(0,\cdot)+ \widetilde{G}_{i+1}^{+,\, k-1}, \quad \mbox{ in } \,  \ie{0,N}.
\end{aligned}
\end{equation}
The remainders are given by
 \begin{equation}\label{eq:errorSWRbis}
\begin{aligned}
&F_i^k =f(U_i^k,\dtms U_i^k,\dxz  U_i^k)-f(\bar{U}_i,\dtms U_i^k,\dxz  \bar{U}_i)+ \bar{U}_i^k  \mbox{ for $n \ge 1$}, \\
& G_i^{-,\, k}= -\frac{\Delta x}{2}f_i^k(0,\cdot)+ g^-(U_i^k(0,\cdot))-g(\bar{U}_i(0,\cdot)),\\
&\widetilde{G}_{i-1}^{-,\, k} = \frac{\Delta x}{2}f_{i-1}^k(J_{i-1},\cdot)+ g^-(\bar{U}_{i-1}^k(J_{i-1},\cdot))-g^-(\bar{U}_{i-1}(J_{i-1},\cdot)),\\
& G_i^{+,\, k}=  -\frac{\Delta x}{2}f_i^k(J_{i},\cdot)+ g^+(\bar{U}_i^k(J_{i},\cdot))-g^+(\bar{U}_i(J_{i},\cdot)),\\
& \widetilde{G}_{i+1}^{+,\, k} =  \frac{\Delta x}{2}f_{i+1}^k(0,\cdot)+ g^+(\bar{U}_{i+1}^k(0,\cdot))-g^+(\bar{U}_{i+1}(0,\cdot)),
\end{aligned}
\end{equation}
and $ G_i^{\pm,\, k}(0)=0$, $\widetilde{G}_i^{\pm,\, k}(0)=0$. For $n=0$, the centered derivative in $x$ are replaced in the expression of $F_i^k$ by a forward or backward derivative.

We define now a discrete energy as follows. We consider  sequences
of the form $V=\{V(j)\}_{0 \le j \le J}$  in $\R^{J+1}$, and we define a bilinear form on $\R ^{J+1}$ by
\begin{equation}  
  a_\Delta(V,W)=\Delta x(\sum_{j=1}^{J}  \dxm(V)(j)\cdot\dxm(W)(j)+\sum_{j=1}^{J-1}V(j)W(j)) .   \end{equation}
%We also define a particular sum denoted by $\sum\mathbb{'}$, which is given in space by
%\[
%\SSx V(j) = \frac{1}{2}V(0) + \sum_{j=1}^{J-1} V(j) + \frac{1}{2}V(J), 
%\] 
%and analogously in time. 
For a mesh function $V$ of time and space, we define
\begin{equation} \label{eq:enerd}   
  \begin{aligned} 
&    E_{K}(V)(n) = \ds \frac{\Delta x}{2}\sum_{j=1}^{J-1}\bigl((\dtm V(j,n))^2+(\dtm V(j,n+1))^2\bigr),\\
&    E_{P}(V)(n) = a_\Delta(V(\cdot,n),V(\cdot,n-1)), \\[3mm]
&    E   = E_{K} +E_{P}. 
  \end{aligned}  
\end{equation}   
The quantity $E_{K}$ is  a discrete kinetic energy. It is less
evident to identify $E_P$ as discrete potential  energy. The following lemma gives a
lower bound for $E$ under a CFL condition, and hence shows that
$E$ is then indeed an energy. The proof is classical (\cite{GHN:2003:OSW}) and is omitted here.
\begin{lemma}\label{lem:EnergyBoundFromBelow}  
  For any  $n \ge 1$, we have
  \begin{equation}\label{energyboundbelow}  
    E(V)(n)\ge\left(1-\frac{\Delta t^2}{\Delta x^2}-\frac{\Delta t^2}{4}\right)  
      E_{K}(V)(n).
  \end{equation} 
  Hence, under the CFL condition 
  \begin{equation}\label{CFL}
     \frac{\Delta t^2}{\Delta x^2}+\frac{\Delta t^2}{4} < 1,
  \end{equation}
  $E$ is bounded from below by an energy.
\end {lemma}     

The following energy estimate is obtained by a discrete integration by parts:
\begin{lemma}\label{lem:enerestim1}
For any $V$ solution of
\begin{align} 
  \left(\dtp\dtm-\dxp\dxm+1\right)V =F ,  \quad  1 \le j \le J -1, \,1 \le n \le N,\label{eq:pbauxint}%\\
%\ds   (\dtp - \frac{\Delta t}{2}\dxp\dxm)V (j,0)= 0, \quad  1 \le j \le J -1,\label{eq:pbauxinit}
  \end{align}
we have for any $n \ge 1$, 
  \begin{equation}\label{energyestimatetimen} 
   \begin{array}{l} 
 \ds E(V)(n)-E(V)(n-1)+
     \ds \frac{\Delta t}{2}
           [(\widetilde{T}^{+}V(0,n))^2+ (\widetilde{T}^{-}V(J,n))^2]\\
\hspace{10mm}   =\ds \frac{\Delta t}{2}
           [(T^{-}V(0,n))^2+ (T^{+}V(J,n))^2]+
   2 \Delta t \Delta x \ds\sum_{j=1}^J F(j,n)\dtz V(j,n),
    \end{array}
  \end{equation} 
and for $n=0$, 
  \begin{equation}\label{energyestimatetime0} 
   \begin{array}{l} 
 \ds E_K(U)(0) + E(U)(0)+ \ds \frac{\Delta t}{4}[(\widetilde{T}^{+}U(0,0))^2+ (\widetilde{T}^{-}U(J+1,0))^2]\\
    \hspace{10mm}= \ds \frac{\Delta t}{4}
           [(T^{-}U(0,0))^2+ (T^{+}U(J,0))^2] 
    +  a_\Delta(P,P)+ \ds 2 \Delta x  \ds\sum_{j=1}^J  Q(j)\dtp U(j,0).
    \end{array}
  \end{equation} 
\end{lemma}

These estimates are obtained by multiplying (\ref{eq:pbauxint}) with $\dtz V(j,n)$ and integrating by parts \cite{GHN:2003:OSW}. \\
We now state  the main result of this section.
\begin{theorem}
Suppose that $f$ is affine with respect to $\partial_x u$. Defining the quantities
\begin{equation}\label{eq:remainder}
     \begin{aligned} 
& R_i^k(n)  =& \widetilde{G}_{i}^{+,\, k}(n)(\widetilde{G}_{i}^{+,\, k}(n)+2 \widetilde{T}^{+}\bar{U}_i^k(0,n))+\widetilde{G}_{i}^{-,\, k}(n)(\widetilde{G}_{1}^{-,\, k}(n)+2\widetilde{T}^{-}\bar{U}_i^k(J_i,n))\\
&           &-G_i^{-,\, k}(n)(G_i^{-,\, k}(n)+2T^{-}\bar{U}_i^k(0,n))- G_i^{+,\, k}(n)( G_i^{+,\, k}(n)+2T^{+}\bar{U}_i^k(J_i,n)),
    \end{aligned}
\end{equation}
and assuming that there exists a positive constant $M$ such that for any iteration number $K$, any domain $\Omega_i$ and any discrete time $n$, the following estimate holds
\begin{equation}\label{eq:hyp}
2  \Delta x \ds\sum_{j=1}^J F_i^k(j,n)\dtz \bar{U}_i^k(j,n)+ \frac{1}{2} R_i^k(n) \le M E(\bar{U}_i^k)(n),
\end{equation} 
then, for $\Delta t$ sufficiently small, the discrete Schwarz algorithm converges in the energy norm.
\end{theorem}
\begin{proof}
We   apply (\ref{energyestimatetimen},\ref{energyestimatetime0}) to $\bar{U}_i^k$. Since the initial data vanish, every term in (\ref{energyestimatetime0}) vanish. Thus $E(\bar{U}_i^k)(0)=0$, and we rewrite \eqref{energyestimatetimen} as
\begin{equation}\label{eq:enestim1}
   \begin{aligned} 
&   \ds  E(\bar{U}_i^k)(n)-E(\bar{U}_i^k)(n-1)  
 \ds +\frac{\Delta t}{2} 
           [(\widetilde{T}^{+}\bar{U}_i^k(0,n)+ \widetilde{G}_{i}^{+,\, k}(n))^2+ (\widetilde{T}^{-}\bar{U}_i^k(J_i,n)+ \widetilde{G}_{i}^{-,\, k}(n))^2]\\
&\hspace{40mm} \ds =\frac{\Delta t}{2}
           [(T^{-}\bar{U}_i^k(0,n) + G_i^{-,\, k}(n))^2+ (T^{+}\bar{U}_i^k(J_i,n)+ G_i^{+,\, k}(n))^2]\\
   &\hspace{40mm} +2 \Delta t \Delta x \ds\sum_{j=1}^{J_i} f_i^k(j,n)\dtz \bar{U}_i^k(j,n)+\frac{\Delta t}{2}R_i^k(n).
    \end{aligned}
\end{equation} 
We now insert the transmission conditions (\ref{eq:errorSWRbc}), translate the indices  in the righthand side, and add the contributions of all subdomains. We define a total internal energy and a total boundary energy as
\[
E_I^k(n)=\sum_{i=1}^IE(\bar{U}_i^k)(n),\ 
E_B^k(n)=\frac{\Delta t}{2} \sum_i [(\widetilde{T}^{+}\bar{U}_i^k(0,n)+ \widetilde{G}_{i}^{+,\, k}(n))^2+ (\widetilde{T}^{-}\bar{U}_i^k(J_i,n)+ \widetilde{G}_{i}^{-,\, k}(n))^2].
\]
With these notations we can write
\begin{equation}\label{eq:enestim3}
   \ds  E_I^k(n)-E_I^k(n-1)+ E_B^k(n)\le E_B^{k-1}(n)+  \ds\sum_{i=1}^I[2 \Delta t \Delta x \ds\sum_{j=1}^{J_i} f_i^k(j,n)\dtz \bar{U}_i^k(j,n)+ \frac{\Delta t}{2} R_i^k(n)].
\end{equation} 
We now sum up (\ref{eq:enestim3}) for $1 \le k \le K$, and define $\hat{E}_I^K(n)= \sum_{k=1}^KE_I^k(n)$:
\begin{equation}\label{eq:enestim4}
   \ds  \hat{E}_I^K(n)-\hat{E}_I^K(n-1)+ E_B^K(n) \le E_B^{0}(n)+\Delta t\sum_{k=1}^K\sum_{i=1}^I[2\Delta x \ds\sum_{j=1}^{J_i} f_i^k(j,n)\dtz \bar{U}_i^k(j,n)+ \frac{1}{2} R_i^k(n)].
\end{equation} 
Under assumption (\ref{eq:hyp}) we deduce that 
\begin{equation}\label{eq:enestim5}
   \ds  \hat{E}_I^K(n)-\hat{E}_I^K(n-1) \le F_{0}(n)+ M \Delta t\hat{E}_I^K(n).
\end{equation} 
The recursive inequality is easy to solve. For $\Delta t$ sufficiently small,  $M \Delta t < 1$, and since $\hat{E}_K(0)=0$, we get
\begin{equation}\label{eq:enestim6}
   \ds  \hat{E}_I^K(n) \le e^{M n\Delta t}\sum_{p=1}^n F_0(p) \le e^{M T}\sum_{p=1}^N F_0(p).
\end{equation} 
This proves that $\sup_{n\le N}\hat{E}_I^K(n)$ is bounded as $K$ tends to infinity. Therefore we have
\begin{equation}\label{eq:enestim7}
\forall n \le N,    \ds \sum_{i=1}^IE(U_i^k-\bar{U})(n) \rightarrow 0 \mbox{ as } k\rightarrow + \infty  
\end{equation} 
which concludes the proof.
\end{proof} 

We are able to prove the assumption (\ref{eq:hyp}) in the case of the linear transmission conditions.  
\begin{corollary}\label{cor:yo}
  Suppose that $f$ is affine with respect to $\partial_x u$. For
  $\Delta t$ sufficiently small, there exists a time $T$ such that the
  discrete Schwarz waveform relaxation algorithm
  (\ref{eq:NumericalSchemeSWR}) with
  linear transmission conditions (\textit{i.e.} $g^\pm=0$) converges to the discrete
  approximation $\bar{U}$ of problem (\ref{clanonlin}) on $(0,T)$.
\end{corollary}

\begin{proof}
Here the remainder $ R_i^k(n)$ reduces to
\begin{equation}\label{eq:remainder0}
     \begin{aligned} 
& R_i^k(n)  =& 2 \Delta x (-f_i^k(0,n)\dtz \bar{U}_i^k(0,n) +f_{i}^k(J_{i},n)\dtz\bar{U}_i^k(J_i,n)),
    \end{aligned}
\end{equation}
and the estimate (\ref{eq:hyp})  amounts to proving that for any $k \le K$, $n \le N$, 
\begin{equation}\label{eq:hyp0}
 \ds\Delta x  \ds\sum_{j=0}^{J_i} f_i^k(j,n)\dtz \bar{U}_i^k(j,n)  \le M_1  E(\bar{U}_i^k)(n).
\end{equation} 
If $f$ is globally Lipschitz in both variables, this is merely an application of the discrete Cauchy-Schwarz lemma. 
\end{proof}

%=================================================================
\section{Numerical Results}
\label{sec:nr}
%=================================================================

%=================================================================
\subsection{Remarks on Overlapping versus Nonoverlapping Schwarz Waveform Relaxation Algorithms and Variants}
\label{subsec:os}
%=================================================================
The original Schwarz algorithm uses overlapping domains (domains $\Omega_i=(a_i,b_i)$ for $1 \le i \le I$, with $a_i < a_{i+1}< b_i< b_{i+1}$), with an exchange of Dirichlet data on the boundary. It is known for elliptic problems to converge, but the smaller the overlap, the slower the convergence \cite{Lions:1988:SAM}. Due to the finite speed of propagation, it converges in a finite number of iterations, given by $N=\lceil cT/L \rceil$ where $c$ is the wave speed, and $L$ the size of the overlap (see \cite{GHN:2003:OSW} for the linear wave equation). Therefore the convergence can be very slow. In the linear case, using absorbing boundary conditions instead of the Dirichlet transmission, even without overlap, improves drastically the convergence, giving in one dimension a number of iterations equal to 2 for some $T$ \cite{GHN:2003:OSW}. In the nonlinear case, we have theoretical convergence results on the linear algorithm for sufficiently small $T$. However, we will see that the latter performs very well on a large time interval, and that the nonlinear algorithm performs even better.\\ 

Our experiments concern the space domain $\Omega=[0,4]$, simulating
$\R$ with linear absorbing boundary conditions at each boundary. The
time interval is $[0,2]$. $\Omega$ is divided in two subdomains. The
initial value is supported in $(0,2)$, with $p(x)=x^3(2-x)^3$, the
initial velocity is $q(x)=3x^2(2-x)^2(x-1)$. This is a good test since
the solution is supported in the first subdomain at $t=0$ and escapes
in the second domain before the end of the computation. 

Note that the Schwarz algorithm can be viewed as  a fixed point algorithm applied to the interface problems 
\[
( H_{1}^{+}, H_{2}^{-})\mapsto (\widetilde{B}^+U_{2}(0,\cdot),\widetilde{B}^-U_{1}(J_1,\cdot))={\cal A}( H_{1}^{+}, H_{2}^{-}).
\]
In all cases, the stopping criterion in the algorithm will be on the
residual $res^k$ for ${\cal A}$.  We also compute the exact discrete
solution in $\Omega$, and measure the discrete global error $E^k$ in
$L^\infty(0,T,L^2(\Omega))$.

\subsection{The Classical Overlapping Schwarz Algorithm}
In this case, $\Omega_1=(0,2+L)$ and $\Omega_2=(2,4)$. The stopping criterion pertains to the residual for the interface problem:  
\[res^k=\biggl(\Vert(u_1^{k+1}-u_1^{k})(2+L,\cdot)\Vert_{L^2}^2+\Vert(u_2^{k+1}-u_2^{k})(2,\cdot)\Vert_{L^2}^2\biggr)^{1/2}.
\] 
We run the computation until the residual $res^k$ is equal to zero. The theoretical minimal number of iterations for the discrete algorithm is $N_{th}=\ds \lceil \frac{\Delta x}{\Delta t}\frac{T}{L}\rceil $. \\

We start with the nonlinear term $f(u,u_t,u_x)=u^3$. Table \ref{table:sc} gives  the number of iterations $N_{comp}$  needed to achieve convergence (\textit{i.e} the error is zero), together with $N_{th}$. 
On the left we  choose the overlap equal to 8 grid points, and vary $\Delta t=1/120$ and $\Delta x=1/100$ to fulfill the CFL condition. On the right we fix   $\Delta t=1/120$ and $\Delta x=1/100$, and vary the overlap $L$.\\

\begin{table}[H]
\begin{minipage}[]{0.45\textwidth}
\centering
\begin{tabular}{|c|ccc|}
\hline
$\Delta x$ & 1/100  &  1/200  &  1/400\\
\hline
$\Delta t$ & 1/120  &  1/240  &  1/480\\
\hline
$N_{comp}$ &29  &    54 &   105  \\
\hline
$N_{th}$ &30  &    60 &  120  \\
\hline
\end{tabular}
\\[2mm]
$L=8\Delta x$
\end{minipage}
\begin{minipage}[]{0.45\textwidth}
\centering
\begin{tabular}{|c|cccc|}
\hline
overlap $L$ & 2 $\Delta x$  & 4$\Delta x$  &  8$\Delta x$  &   16$\Delta x$\\
\hline
$N_{comp}$ &108  &    55 &  28 &  15  \\
\hline
$N_{th}$    &121  &    61 &   31 &    16\\
\hline
\end{tabular}
 \\[2mm]
$\Delta x=1/100$, $\Delta t=1/120$ 
\end{minipage}\\
\caption{Number of iterations to achieve convergence for the classical  Schwarz algorithm}
\label{table:sc}
\end{table}

\begin{figure}[H]
  \centering
  \includegraphics[width=0.45\textwidth]{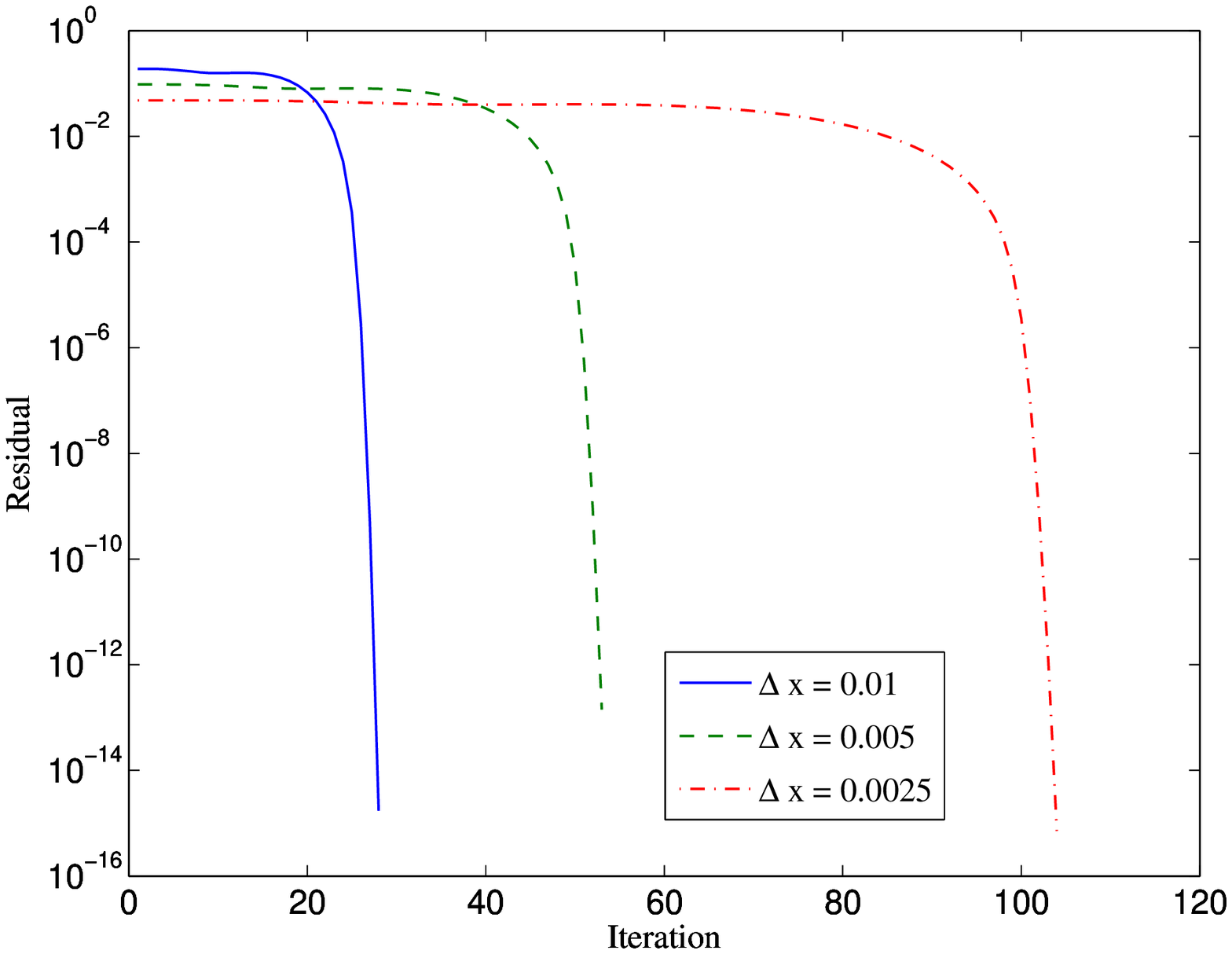}
  \includegraphics[width=0.45\textwidth]{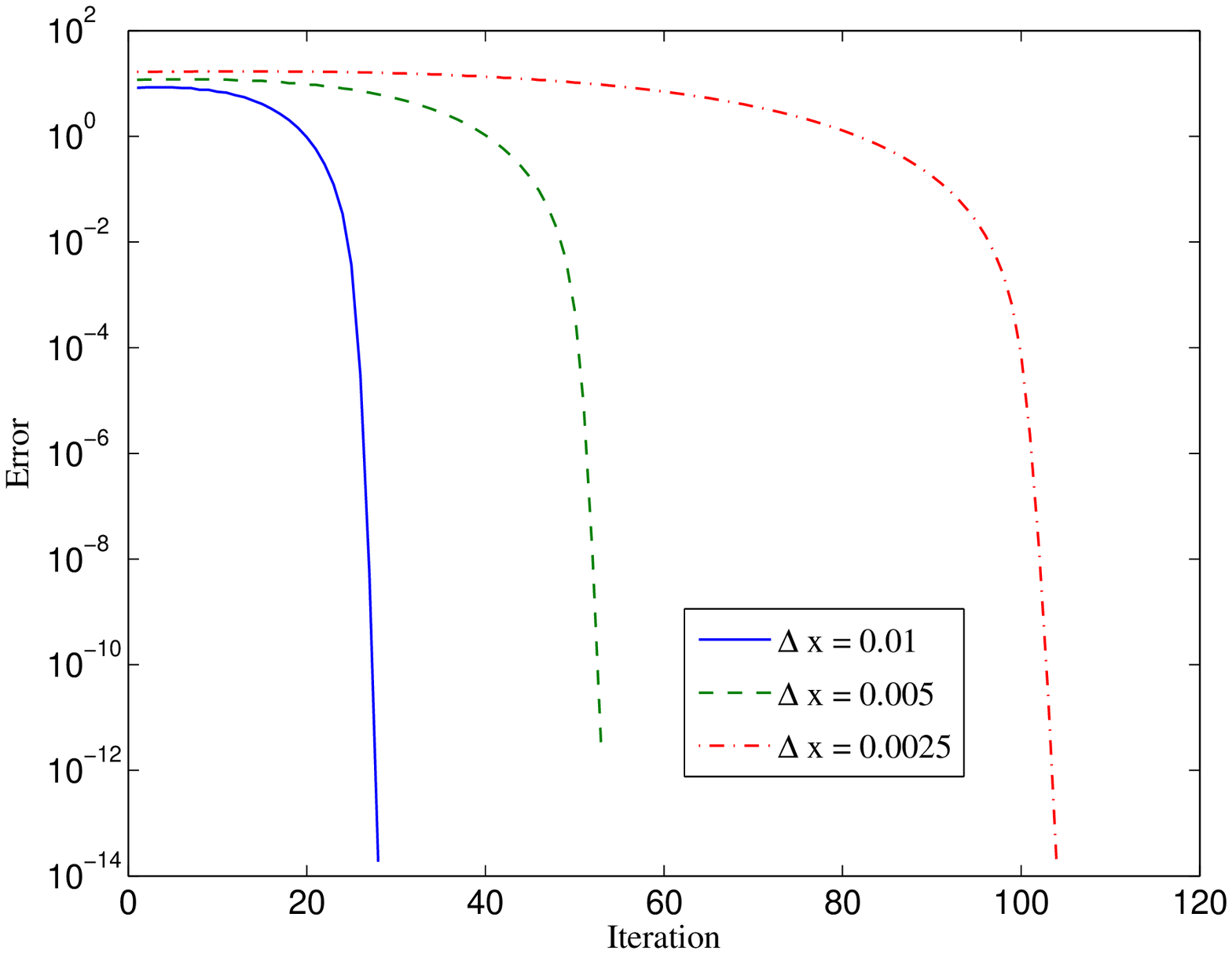}
\caption{ Variation of the residual  $res^k$ (left) and the error $E^k$(right) as a function of the iteration number $k$. }
\label{fig:sc1}
 \end{figure}

We show on Figure \ref{fig:sc1} the convergence history of the classical Schwarz algorithm for various values of the mesh size and an overlap equal to eight gridpoints.We check in each case that the error $E^k$ vanishes together with the residual. Furthermore we can see that the error decays very slowly for many iterations, and reaches zero in a few iterations, independently of the mesh size (three or four in all cases). The behaviour is very similar to what happens for the linear wave equation : only the finite speed of propagation produces convergence, which takes place when the signal has left the domain.  The algorithm behaves similarly for other nonlinearities.

 \subsection{The Non linear Nonoverlapping Schwarz Algorithms}
We will see in this section that our strategy greatly improves the
performances of the Classical Schwarz algorithm. Note that, whereas the classical algorithm only
converges in the presence of an overlap, and the smaller the overlap, the slowlier the convergence, our algorithms are run without overlap, 

Here the residual is given by:  
\[
res^k=\biggl(\Vert(G^{+,\, k+1}-G^{+,\, k})(a_1,\cdot)\Vert_{L^2}^2+\Vert(G^{-,\, k+1}-G_2^{-,\, k})(a_1,\cdot)\Vert_{L^2}^2\biggr)^{1/2}.
\] 
We start with $f=u^3$. In this case the transmission operators are the
same and linear, since $g^{\pm}=0$. We have proved in Corollary
\ref{cor:yo} that there exists a final time $T$ for which the discrete
algorithm is convergent.

\begin{figure}[H]
  \centering
\begin{tabular}[c]{ccc}
  \includegraphics[width=0.30\textwidth]{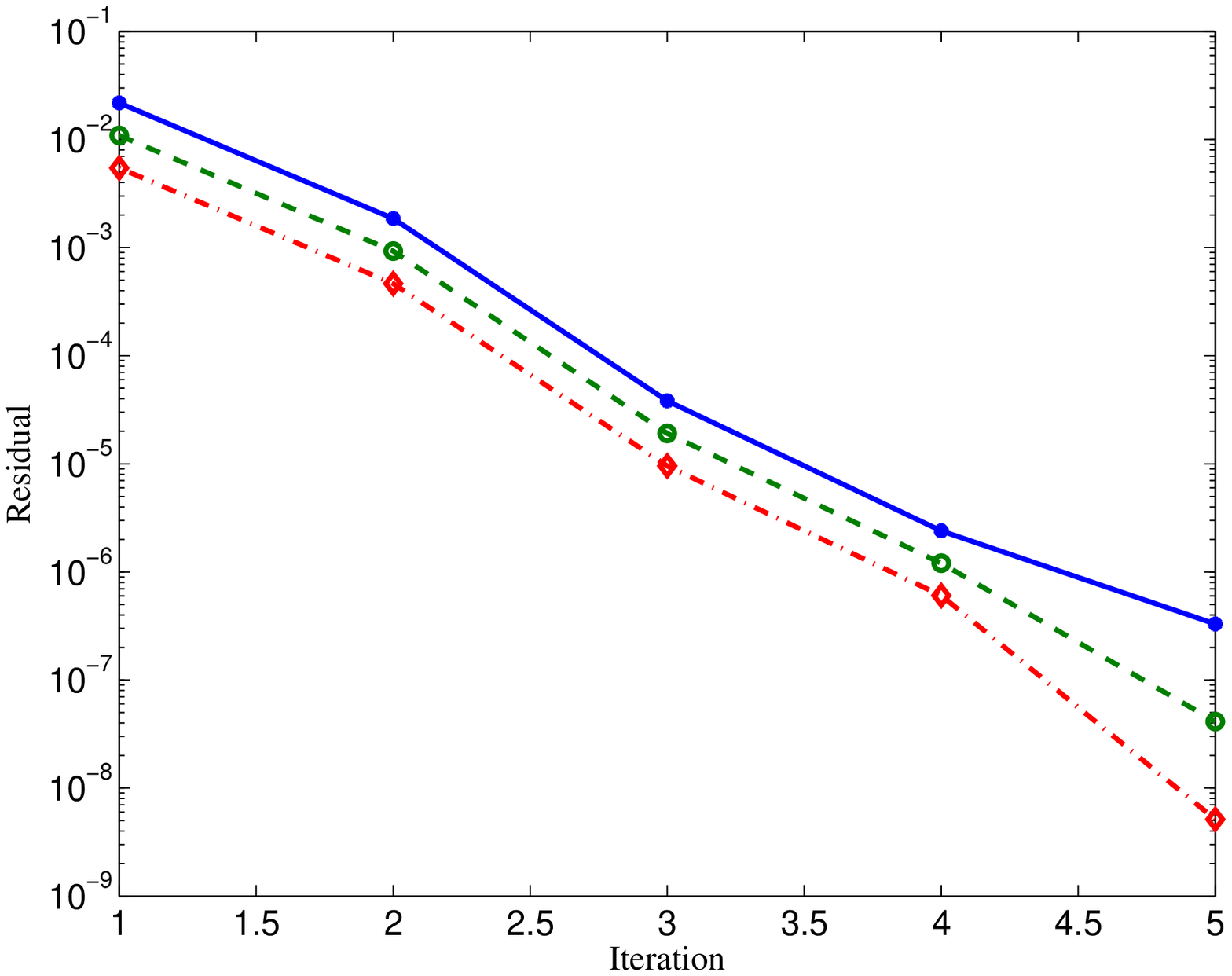}&
   \includegraphics[width=0.30\textwidth]{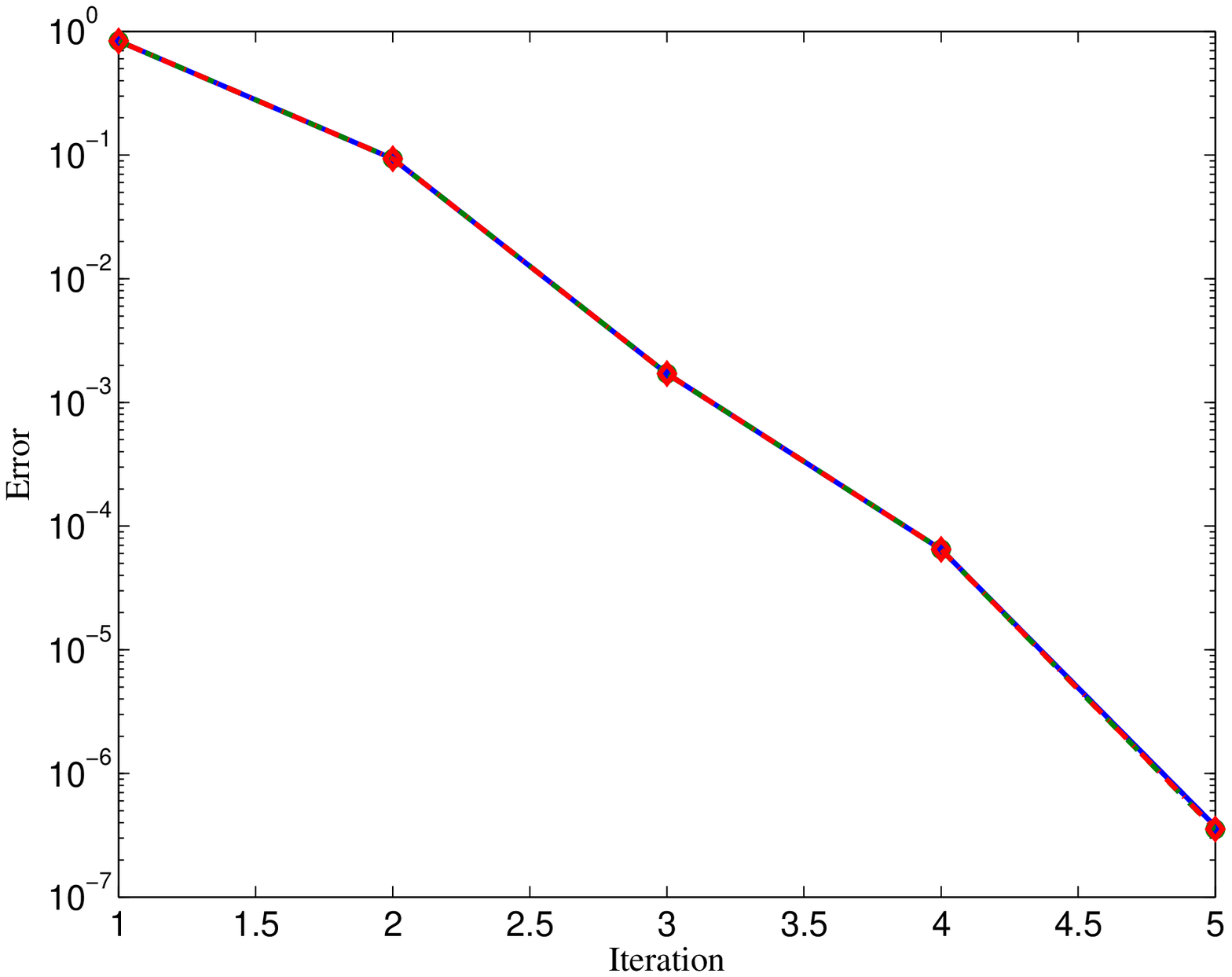}
\end{tabular}
 \caption{$f=u^3$. Left residual, right: error. $*$:  $\Delta x=1/100$, $\Delta t=1/120$, o:  $\Delta x=1/200$, $\Delta t=1/240$, $\diamondsuit$:  $\Delta x=1/200$, $\Delta t=1/240$.}
\label{fig:snl3201short}
 \end{figure}

In the forthcoming computations, the theoretical and numerical data
are the same as before.  Figure \ref{fig:snl3201short} plots the
convergence history for various mesh sizes. The computation is stopped
as soon as the residual reaches $0.5 \,10^{-7}$.
\noindent We see that the algorithm converges very rapidly, independently of the mesh size.  

\begin{figure}[H]
  \centering
\begin{tabular}{ccc}
  \includegraphics[width=0.30\textwidth]{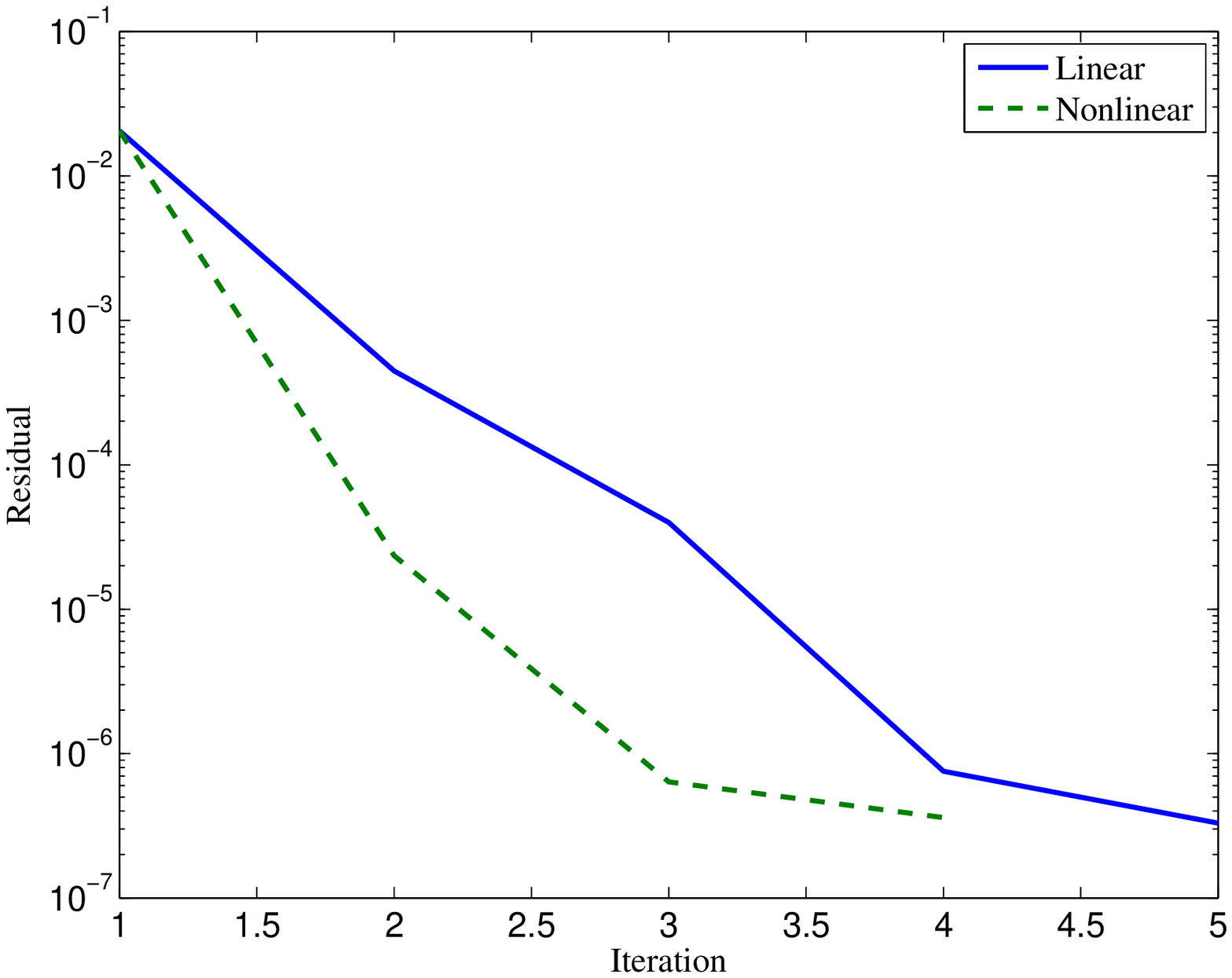}&
  \includegraphics[width=0.30\textwidth]{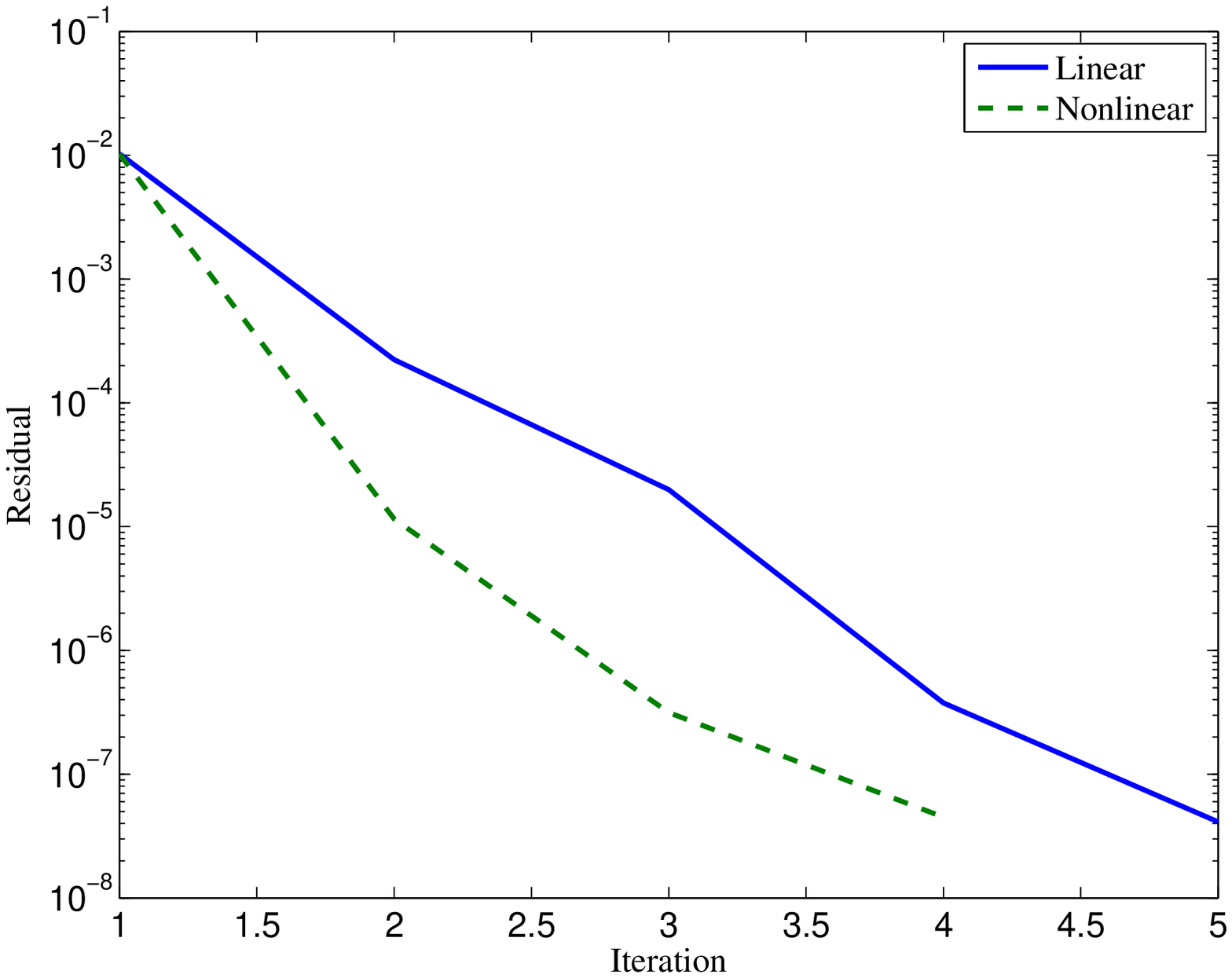}&
  \includegraphics[width=0.30\textwidth]{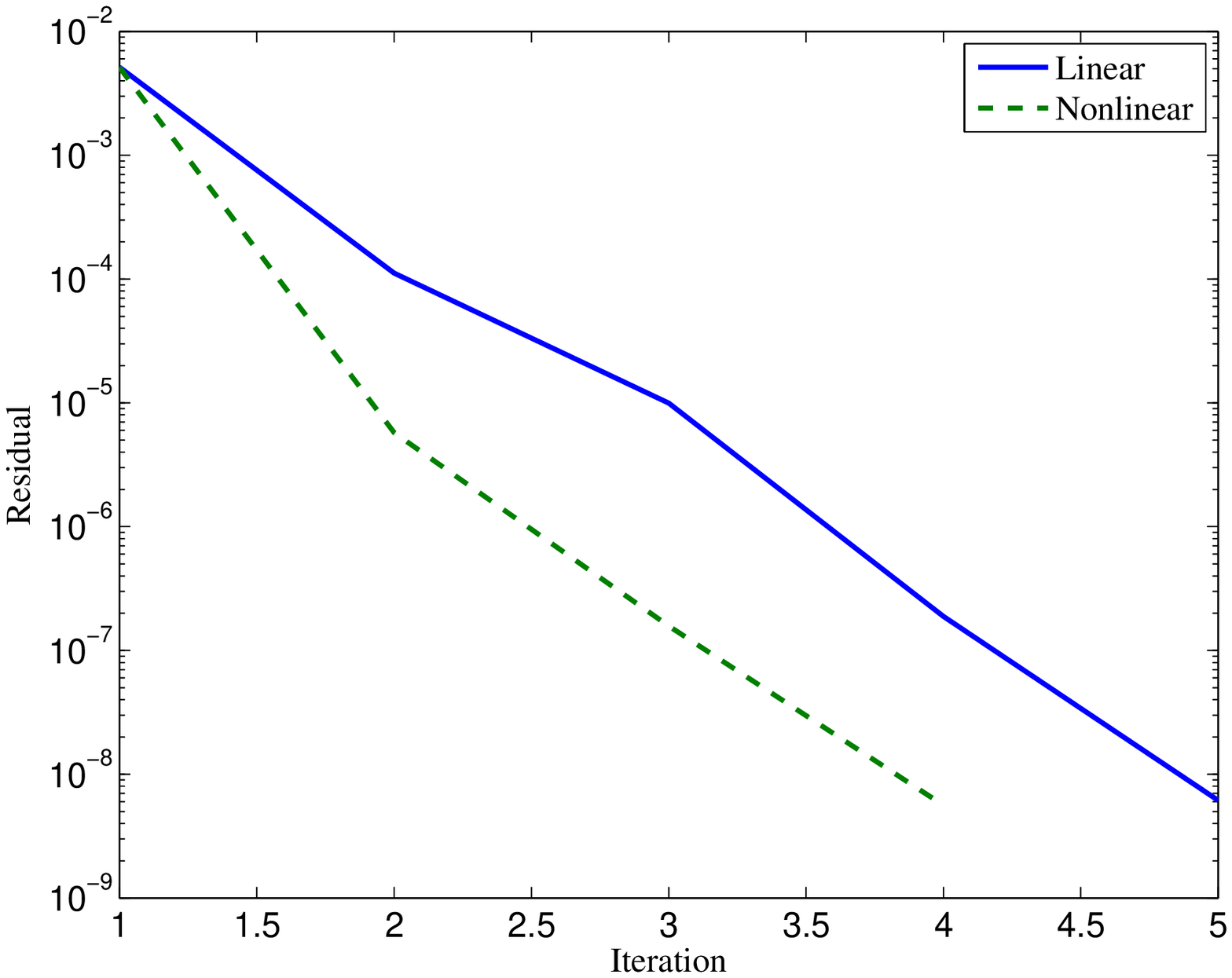}\\
  \includegraphics[width=0.30\textwidth]{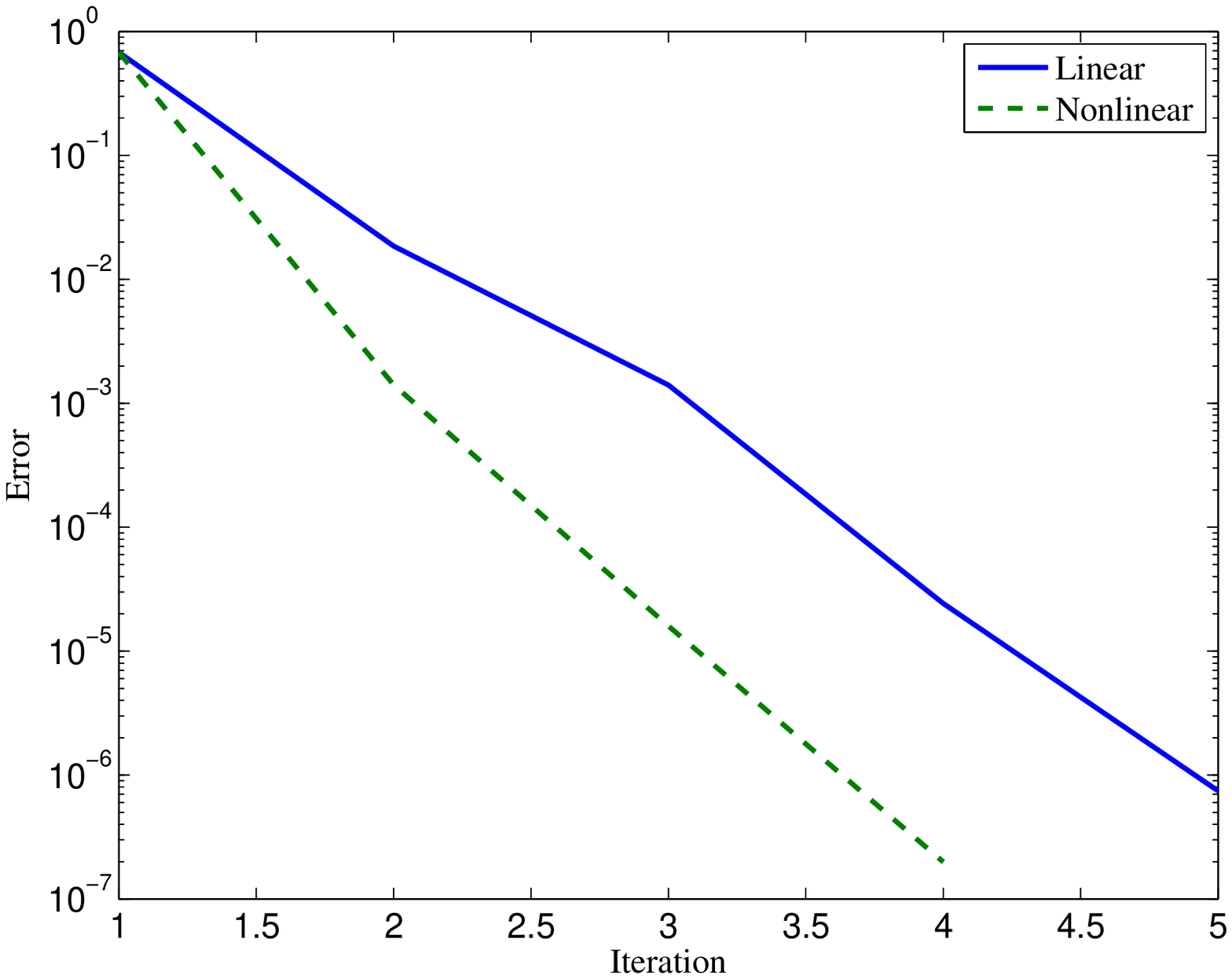}&
  \includegraphics[width=0.30\textwidth]{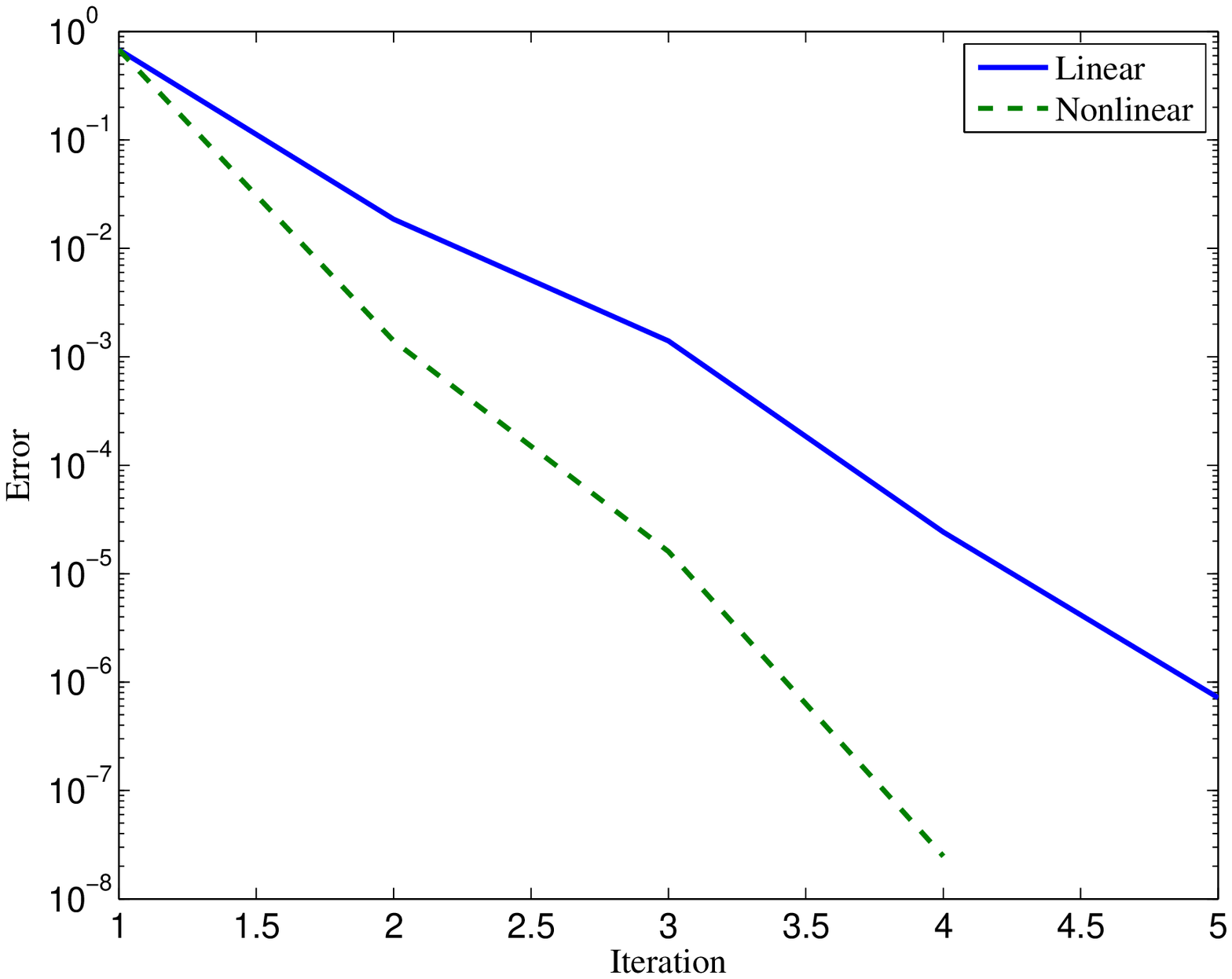}&
  \includegraphics[width=0.30\textwidth]{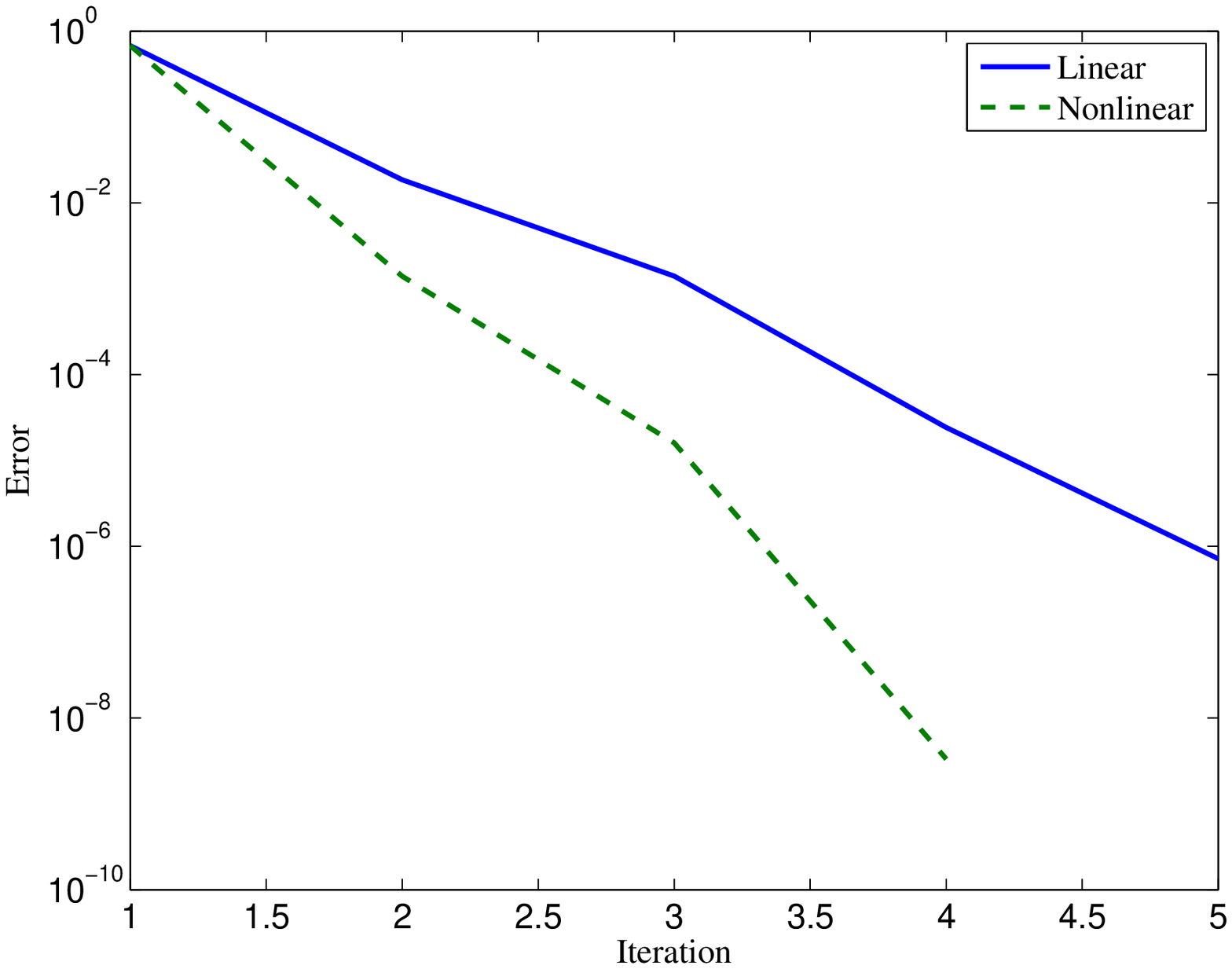}\\
 $\Delta x=1/100$, $\Delta t=1/120$&
 $\Delta x=1/200$, $\Delta t=1/240$&
 $\Delta x=1/400$, $\Delta t=1/480$
\end{tabular}
 \caption{$f=u^2u_x$. Top: residual, bottom: error. Solid: linear, dash: nonlinear.}
\label{fig:snlu2uxs}
 \end{figure}

We consider now the case $f=u^2u_x$, and nonlinear transmission
conditions, \textit{i.e.} $g^{\pm}(u)=\pm u^3/6$.  In Figure
\ref{fig:snlu2uxs}, we plot the convergence history for various mesh
sizes. The computation is stopped as soon as the residual reaches $0.5
\,10^{-7}$.

\noindent Both the linear and the nonlinear transmission conditions behaves very well. 
The convergence takes
place in 5 iterations with the linear transmission condition, in 4
iterations with the nonlinear one. 

\begin{figure}[H]
  \centering
  \includegraphics[width=0.30\textwidth]{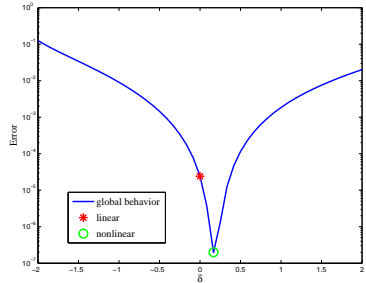}
 \caption{$f=u^2u_x$. Variations of the error as a function of the nonlinearity coefficient $\delta$}
\label{fig:comp}
 \end{figure}

We also can vary the nonlinearities in the transmission conditions, we use a real parameter $\delta$, and 
$g^{\pm}=\pm \delta u^3$. The nonlinear strategy
corresponds to $\delta=1/6$, whereas the linear one is obtained for
$\delta=0$. We draw in Figure \ref{fig:comp} the error curves after 3
iterations for the same initial values as before, the mesh sizes are
$\Delta x=1/100$, $\Delta t=1/120$. We observe that the nonlinear
strategy corresponds precisely to the optimal numerical value of the parameter $\delta$, validating the high frequency approach.
We have carried out the same computations in the case $f=u^2u_t$. We do not display the results here since they are very similar.

\section{Conclusion}
We have presented a linear and a nonlinear Schwarz waveform relaxation
algorithm without overlap for the semilinear wave equation. On the
continuous level, we proved the convergence for sufficiently small
time intervals. We designed a discrete algorithm, in such a way that, if convergent, the algorithm converges to the discrete solution in the whole domain, which we prove when the nonlinearity is affine in $\partial_x u$. In that case we proved the convergence for the linear transmission condition. Numerical
experiments highlight the fast convergence to the discrete full domain
solution in a large time domain in both linear and nonlinear
strategies, without overlap. Furthermore, we have shown that our
nonlinear transmission conditions give optimal results within a large
class of transmission conditions.

\section*{Acknowledgments}
The authors are very grateful to Pr. Martin Gander from Gen\`eve
University whose matlab scripts were a basis for the present work.
\appendix

\bibliographystyle{plain} 
\bibliography{ddnl}
%\bibliography{ddm}

\end{document}